\documentclass[12 pt]{article}
\usepackage{color}
\usepackage{tcolorbox}
\usepackage{setspace}
\usepackage{mathtools}
\usepackage{blindtext}
\usepackage{float}
\usepackage{parskip}
\usepackage{bbold}
\usepackage{mathtools, amssymb}
\usepackage{tikz}
\usepackage{arxiv}
\usepackage[ruled,longend]{algorithm2e}
\usepackage[utf8]{inputenc} 
\usepackage[T1]{fontenc}  
\usepackage[linkbordercolor=white]{hyperref}
\usepackage{url}       
\usepackage{booktabs}      
\usepackage{amsfonts}       
\usepackage{nicefrac}       
\usepackage{microtype}     
\usepackage{titlesec}
\usepackage{graphicx}
\usepackage{changepage}
\usepackage[shortlabels]{enumitem}
\usepackage{doi}
\usepackage{amsmath, nccmath}
\usepackage[ruled,longend]{algorithm2e}
\graphicspath{ {images/} }

\title{Explicit solution to an optimal two-player switching game in infinite horizon}

\author{ {\includegraphics[scale=0.06]{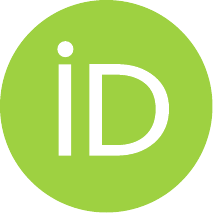}\hspace{1mm}Brahim EL ASRI }\\
	Equipe Aide à la Décision,\\
	 Université Ibn Zohr, \\
	ENSA, BP 1136, Maroc, \\
	\texttt{b.elasri@uiz.ac.ma} \\
	\And
	{\includegraphics[scale=0.06]{orcid.pdf}\hspace{1mm} Magnoudéwa PAKA } \\
	Equipe Aide à la Décision,\\
	 Université Ibn Zohr, \\
	ENSA, BP 1136, Maroc, \\
	\texttt{magnoudewa.paka@edu.uiz.ac.ma} \\
}

\renewcommand{\headeright}{}
\renewcommand{\undertitle}{}
\renewcommand{\shorttitle}{Explicit solution to an optimal two-player switching game in infinite horizon}
\hypersetup{
pdftitle={A template for the arxiv style},
pdfsubject={q-bio.NC, q-bio.QM},
pdfauthor={David S.~Hippocampus, Elias D.~Striatum},
pdfkeywords={First keyword, Second keyword, More},
}
\begin{document}
\hypersetup{pdfborder=0 0 0}
\maketitle
\newtheorem{theo}{Theorem}[section]
\newtheorem{problem}{Problem}[section]
\newtheorem{pro}{Proposition}[section]
\newtheorem{cor}{Corollary}[section]
\newtheorem{axiom}{Definition}[section]
\newtheorem{rem}{Remark}[section]
\newtheorem{lem}{Lemma}[section]
\newtheorem{ex}{Example}[section]
\newtheorem{proof}{Proof}[section]
\newtheorem{Heuristics}{Heuristics}
\newtheorem{ass}{Assumption}[section]
\newcommand{\eqskip}{ \qquad \qquad \qquad \qquad}

\newcommand{\bass}{\begin{ass}}
\newcommand{\eass}{\end{ass}}
\newcommand{\bpf}{\begin{proof}}
\newcommand{\epf}{\end{proof}}
\newcommand{\brm}{\begin{rem}}
\newcommand{\erm}{\end{rem}}
\newcommand{\bethe}{\begin{theo}}
\newcommand{\eethe}{\end{theo}}
\newcommand{\bl}{\begin{lem}}
\newcommand{\el}{\end{lem}}
\newcommand{\bp}{\begin{pro}}
\newcommand{\ep}{\end{pro}}
\newcommand{\bcor}{\begin{cor}}
\newcommand{\ecor}{\end{cor}}
\newcommand{\be}{\begin{equation}}
\newcommand{\ee}{\end{equation}}
\newcommand{\beq}{\begin{eqnarray*}}
\newcommand{\eeq}{\end{eqnarray*}}
\newcommand{\beqa}{\begin{eqnarray}}
\newcommand{\eeqa}{\end{eqnarray}}
\newcommand{\dg}{\displaystyle \delta}
\newcommand{\cm}{\cal M}
\newcommand{\cF}{{\cal F}}
\newcommand{\cR}{{\cal R}}
\newcommand{\bF}{{\bf F}}
\newcommand{\tg}{\displaystyle \theta}
\newcommand{\w}{\displaystyle \omega}
\newcommand{\W}{\displaystyle \Omega}
\newcommand{\vp}{\displaystyle \varphi}
\newcommand{\ig}[2]{\displaystyle \int_{#1}^{#2}}
\newcommand{\integ}[2]{\displaystyle \int_{#1}^{#2}}
\newcommand{\produit}[2]{\displaystyle \prod_{#1}^{#2}}
\newcommand{\somme}[2]{\displaystyle \sum_{#1}^{#2}}
\newlength{\inter}
\setlength{\inter}{\baselineskip}
\setlength\parindent{24pt}
\setlength{\baselineskip}{7mm}
\newcommand{\no}{\noindent}
\newcommand{\rw}{\rightarrow}
\def \ind{1\!\!1}
\def \R{I\!\!R}
\def \N{I\!\!N}
\def \cadlag {{c\`adl\`ag}~}
\def \esssup {\mbox{ess sup}}
\titlespacing*{\section}{0pt}{0.3\baselineskip}{0.5\baselineskip}
\titlespacing*{\subsection}{0pt}{0.3\baselineskip}{0.5\baselineskip}
 \renewcommand{\baselinestretch}{0.7}
{\setstretch{1.4}
\begin{abstract}
\qquad In this paper we use viscosity approach to provide an explicit solution to the problem of a two - player switching game.
We characterize the switching regions which  reduce the switching problem into one of finding a finite number of threshold values in
state process that would trigger switchings and then derive  an explicit solution to this problem. 
The state process is a one dimensional Itô diffusion process and switching costs are allowed to be non-positive.
We also suggest a numerical procedure to compute the value function in case we know the qualitative structure of switching regions
 and we illustrate our results by numerical simulations.
\end{abstract}
\medskip
\\\textbf{AMS Classification subjects}: 60G40 ; 62P20 ; 91B99 ; 91B28 ; 35B37 ; 49L25.
\medskip

\section{Introduction}
\setcounter{equation}{0}
\renewcommand{\theequation}{1.\arabic{equation}}
\no  Differential games are problems in which we model and analyze a conflict in the context of dynamical system. Such problems were addressed in many fields including economics, finance, energy and can be addressed in any other field where there are individuals subject to decisions according to a state variable that evolves according to a differential equation. A typical example in finance explored by K. Suzuki \textcolor{blue}{\cite{K2}} is the case of a pair trading strategy in a context of mean reverting asset portfolio: Consider two similar stocks
that are highly correlated. Assume that the spread between the two
stock prices fluctuates randomly and the spread has a long-run
mean. Sometimes the spread process diverges from the long-run
mean and sometimes it converges. If the spread widens, the expensive stock is sold and the cheap stock is purchased. As the spread narrows again, profit is taken by unwinding the pairs position.
A typical example in economics is the differential game of capitalism explored by K. Lancaster \textcolor{blue}{\cite{L}} in which we assume two players:  the workers and the capitalists. The workers control the share of their consumption in total output while the capitalists control the share of investment in the surplus. In energy, we can think of a manager of a power plant which produces electricity and which
can operate in  several modes of production. The price of inputs in each mode fluctuates and switching from one regime to another incurs cost. The strategy of the manager is to select a sequence of switching times and the associated mode so as to maximize her profit.

 In this paper, we study a two player switching game. To fix ideas, we consider the state process of the stochastic differential game, defined as the solution of the following stochastic equation:
\begin{equation}
\label{az}
\begin{array}{ll}
 X_{s}=x+\integ{0}{s}b(r,X_{r},\mu_r,\nu_r)dr+\integ{0}{s}\sigma(r,X_{r},\mu_r,\nu_r)dW_{r}\qquad s\geq t,
\end{array}
\end{equation}
with $X_{0^-} = x$ the initial state. Here $W$ is a d-dimensional Wiener process, while
\begin{equation*}
\mu_t=\sum_{n>0} \xi_{n} \mathbb{1}_{\tau_n \leq t < \tau_{n+1}} \quad \quad \text{and} \quad \quad \nu_t=\sum_{m>0} \eta_{m} \mathbb{1}_{\rho_m \leq t < \rho_{m+1}},
\end{equation*}
represent the regime of  player I (resp Player II)  at each time, $(\tau_n)_{n>0}$ and $(\rho_m)_{m>0}$, their respective sequence of switching times, $\xi_{n}$ the regime of player I just after switching at $\tau_n$ and $\eta_{m}$, the regime of player II just after switching at $\rho_m$. The payoff to be controlled is:
\begin{equation}
\label{eq00}
\begin{array}{ll}
\mathbb{E}\bigg[\integ{0}{\infty}e^{-rs}f(X^{\mu_s,\nu_s}_{s},\mu_s,\nu_s)ds-
\sum\limits_{n\geq 1}e^{-r\tau_n} C(\xi_{n-1},\xi_n) +\sum\limits_{m\geq 1}e^{-r\rho_m}\chi(\eta_{m-1},\eta_m)\bigg].
\end{array}
\end{equation}
The first player chooses the control $\mu$ from a given finite set of regimes $\mathcal{D}$ to maximize the payoff \textcolor{blue}{(\ref{eq00})} and each of his action is related with one cost $C$, while the second player chooses the control $\nu$ from  set of regimes $\mathcal{D}$ to minimize the payoff \textcolor{blue}{(\ref{eq00})}, and each of his actions is associated with the other cost $\chi$. The 
problem we will investigate is to characterize switching regions for each player and then derive an explicit solution to the game using viscosity approach together with the associated quasi-variational inequalities.
The Isaacs system of
quasi-variational inequalities for this switching game is the following:
for any $i,j\in \mathcal{D}$, and $x\in \mathbb{R}$,
\begin{equation} \label{eq:HJBI0}
\begin{array}{c}
max\Big\{ min\Big[rV_{ij}(x)-\mathcal{L}_{ij}V_{ij}(x)-f(x,i,j);\qquad\qquad\qquad\qquad\\ \qquad V_{ij}(x)-M_{ij}[V](x)\Big];V_{ij}(x)-N_{ij}[V](x)\Big\}=0,
\end{array}
\end{equation}
and
\begin{equation} \label{eq:HJBI01}
\begin{array}{c}
min\Big\{ max\Big[rV_{ij}(x)-\mathcal{L}_{ij}V_{ij}(x)-f(x,i,j);\qquad\qquad\qquad\qquad\\ \qquad V_{ij}(x)-N_{ij}[V](x)\Big];V_{ij}(x)-M_{ij}[V](x)\Big\}=0.
\end{array}
\end{equation}
where,
$$M_{ij}[V](x)= \max_{k\ne i}\{V_{kj}(x)-C(i,k)\}, \qquad N_{ij}[V](x)= \min_{l\ne j}\{V_{il}(x)+\chi(j,l)\}.$$
and
$$\mathcal{L}_{ij}v_{ij}(x)=\frac{1}{2} \mathrm{Tr}[\sigma_{ij} (\sigma_{ij})^* \nabla_x^2 v_{ij}(x)]+ \big < b_{ij} ,\nabla_x v_{ij}(x)\big > .$$

  As part of applied mathematics, differential games subject arouses more and more interest among researchers.  Several authors have conducted their research on this subject, including W. H. Fleming and P. E.  Souganidis \textcolor{blue}{\cite{FS}} who studied the existence of value functions of two-player zero-sum stochastic differential games, C. Evans and E. Souganidis \textcolor{blue}{\cite{ES}} who tackled the  representation formulas for solutions of Hamilton-Jacobi-Isaacs (HJI for short) equations. E. Bayraktar and M. Egami \textcolor{blue}{\cite{EM}} on their side solved explicitly a one-dimensional optimal switching problem by means of dynamic programming principle and the exercise characterisation of the value function. However, the study of two person zero-sum differential games dates back to the work of R. Isaacs \textcolor{blue}{\cite{I}}, and one of the first problem analyzed was the homicidal chauffeur game. In this problem, Isaacs studied  the pursuit problem which pits a hypothetical runner, who can only move slowly, but is highly maneuverable, against the driver of a motor vehicle, which is much faster but far less maneuverable, who is attempting to run him down.  Both runner and driver are assumed to never tire.  The question to be solved is: under what circumstances, and with what strategy, can the driver of the car guarantee that he can always catch the pedestrian, or the pedestrian guarantee that he can indefinitely elude the car. In the finite and infinite horizon framework, the switching game have been studied by several authors. The
most recent works discussing this topic include the papers by Djehiche et al \textcolor{blue}{\cite{DHH}}, Tang and Hou \textcolor{blue}{\cite{THH}}  and B. El Asri, S. Mazid \textcolor{blue}{\cite{BM}}. They showed existence and uniqueness of a continuous viscosity solution of
value function. The explicit solution for switching game has been addressed by Vathana Ly Vath  and Huyen Pham \textcolor{blue}{\cite{LP}} in one player case. 

However the complexity of differential game problems has led researchers to lean towards numerical solutions, so that in the literature, little work is devoted to analytical or explicit solutions. Explicit solutions, although complicated to investigate, have the advantage of being precise, i.e. do not depend on a certain number of iterations to be considered reliable, constraints frequently encountered in numerical  solutions. This fact is  the motivation of our work. The novelty of this paper is to provide an explicit characterization of the switching regions showing when and
where it is optimal to change the regime and derive an explicit solution to thedifferential game problem. To tackle this
problem, we use viscosity approach together with quasi-variational inequalities.
Additionaly, we provide a numerical procedure to compute threshold values in case we know qualitative
structures of switching regions without knowing explicitly the value of the threshold values and we give
some numerical examples with graphic illustrations.

The paper is organised as follow: In section 2, we formulate the problem and state assumptions of the problem with some preliminary results; in section 3, we state the quasi-variational inequalities and switching regions. Section 4 is devoted to the explicit solution of the problem and  finally,in section 5 We suggest a numerical procedure to compute the value function in case we know the qualitative structure of switching regions we illustrate our result with some numerical simulations in section 6.
\section{Problem formulation, Assumptions and Preliminary results}
\setcounter{equation}{0}
\renewcommand{\theequation}{2.\arabic{equation}}

\subsection{Problem formulation}

\no The problem consists of a two-player game: the first player ( player I) and the second player ( player II), both subject to decision making according to the value taken by a state process $X:=(X_t)_{t \geq 0}$. Each of the players has several modes or regimes in which he can switch. Let us denote by $\mathcal{D}:=\{1,2,....,d\}$ the finite set of regimes, i.e., each regime is represented by an index in $\{1,2,...,d \}$. The indicator of the regimes for both players is modeled by a two-dimensional cadlag process $\textit{I}_t$ valued in $\mathcal{D} \times \mathcal{D}$. For each value of $\textit{I}_t$, there is a running profit function $f_{\textit{I}_t}$ that is assumed to be non-negative.  The state process $X$ can be seen as the price of a commodity and therefore is valued in $\mathbb{R}_+^*$. We denote by $J$ the payoff of the game. The player I aims at maximizing  the payoff $J$ by choosing  suitable switching regimes $\xi:=(\xi_n)_{n>0}$ at suitable stopping times $\tau:=(\tau_n)_{n>0}$  known as the optimal sequence  and the player II aims at minimizing the payoff by choosing suitable regimes $\eta:=(\eta_m)_{m>0}$ at suitable stopping times  $\rho:=(\rho_m)_{m>0}$. Switching from a regime $i$ to another regime $j$ incurs a cost denoted by $c_{ij}$ for player I and $\chi_{ij}$ for player II.    The payoff of the game is expressed as follow:
\begin{equation}
\label{eq0}
\begin{array}{rr}
J(x,\xi,\eta) =\mathbb{E}\bigg[\integ{0}{\infty}e^{-rs}f(X^{\mu_s,\nu_s}_{s},\mu_s,\nu_s)ds-
\sum\limits_{n\geq 1}e^{-r\tau_n} C(\xi_{n-1},\xi_n) +\sum\limits_{m\geq 1}e^{-r\rho_m}\chi(\eta_{m-1},\eta_m)\bigg], 
\end{array}
\end{equation}
where $r$ is a positive discount rate,  $\mu_t$ (resp $\nu_t$) is the regime of player I( resp player II) at time $t$ and given by:
\begin{equation*}
\mu_t=\sum_{n>0} \xi_{n} \mathbb{1}_{\tau_n \leq t < \tau_{n+1}} \quad \quad \text{and} \quad \quad \nu_t=\sum_{m>0} \eta_{m} \mathbb{1}_{\rho_m \leq t < \rho_{m+1}}.
\end{equation*}
We can notice that by definition we have $\textit{I}_t= (\mu_t, \nu_t).$  \\
\qquad In this work, the process $X_t$ is of geometric brownian motion type and is expressed as follows:
\begin{equation}
\begin{cases}
dX_t  &=  b(X_t, \mu_t, \nu_t) dt +  \sigma(X_t, \mu_t, \nu_t) dW_t , \\
 X_{0} &=  x, \quad \textit{I}_{0^-}=(i,j), 
\end{cases}
\end{equation}
where $$b(X_t, \mu_t, \nu_t) = b_{\textit{I}_t} X_t \quad \text{ and } \quad \sigma(X_t, \mu_t, \nu_t)=\sigma_{\textit{I}_t}X_t.$$

Throughout this paper, diffusion operators refer to $b_{\textit{I}_t}$ and $\sigma_{\textit{I}_t}$. We shall also denote the solution of the stochastic differential equation of the process as $X^x_{I_0}$ or simply by $X^x$.
\subsection{Assumptions}
Throughout this paper $d$ is an integer.  Let $\mathcal{D}=\{1,...,d\}$ be the finite set of regimes, and assume the following assumptions:  \\
\no\textbf{[H1]}  $b:\mathbb{R}^{+}\times\mathcal{D}\times\mathcal{D}\rightarrow
\mathbb{R}^{+}$ and $\sigma :\mathbb{R}^{+}\times\mathcal{D}\times\mathcal{D}\rightarrow
\mathbb{R}^{+}$ be two continuous functions for which there
exists a constant $C>0$ such that for any
$x,x^{\prime }\in \mathbb{R}^{+}$ and $i,j\in\mathcal{D}$:
\begin{equation}
|\sigma(x,i,j)-\sigma(x^{\prime },i,j)|+|b(x,i,j)-b(x^{\prime},i,j)|\leq
C|x-x^{\prime }|. 
 \label{eqs}
\end{equation}%
Thus they are also of linear growth,  i.e., there exists a constant $C$ such that for any $x\in\mathbb{R}^+$ and $i,j\in\mathcal{D}$:
\begin{equation}
|\sigma(x,i,j)|+|b(x,i,j)|\leq C(1+|x|). 
\end{equation}
\textbf{[H2]}\; $f:\mathbb{R}^{+}\times\mathcal{D}\times\mathcal{D} \rightarrow \mathbb{R}^{+}$
is a continuous function for which there exists a constant $C$ such that for each $i,j\in\mathcal{D}$, $x,x^{\prime }\in \mathbb{R}^{+}$:
 \begin{equation}
 \qquad |f(x,i,j)|\leq C(1+|x|) \qquad\text{and}
\qquad |f(x,i,j)-f(x^{\prime},i,j)|\leq
C|x-x^{\prime }|. \qquad\qquad\quad
\end{equation}
\textbf{[H3]}\; For any $i,j\in\mathcal{D}$, the switching costs $C(i,j)$ and $\chi(i,j)$ are constants, and we assume the triangular condition :
\begin{equation}\label{co1}
 c(i,k)< c(i,j)+c(j,k), \quad j\ne i,k. \quad \text{and } \quad
\chi(i,k)< \chi(i,j)+\chi(j,k), \quad j\ne i,k,
\end{equation}
which means that it is less expensive to switch directly in one step from regime $i$ to $k$ than
in two steps via an intermediate regime $j$.  Notice that a switching costs $c(i,j)$ and $\chi(i,j)$ may be negative,
and conditions \textcolor{blue}{(\ref{co1})}  for $i=k$ prevents an arbitrage by simply switching back
and forth, i.e. 
\begin{equation}\label{co3}
0< c(i,j)+c(j,i), \qquad \text{and}
\qquad
0< \chi(i,j)+\chi(j,i). 
\end{equation}
We will use alternatively $c(i,j)$ or $c_{ij}$ to represent the switching costs and $f_{ij}(. )$ or $=f(. ,i,j)$ represent profit functions.\\
\textbf{[H4]} The family of random variables $(c_n)_{n\geq 1}$ (resp. $(\chi_n)_{n\geq 1}$), where $(c_n)_{n\geq 1}$ (resp. $(\chi_n)_{n\geq 1},$) is the total cost of the first $n$ switches
      $$c_n=\sum\limits_{m= 1}^{n}e^{-r\tau_m} c(\xi_{m-1},\xi_m) \quad\text{(resp}.\; \chi_n=\sum\limits_{l=1}^{n}e^{-r\rho_l}\chi(\eta_{l-1},\eta_l))$$
     converges $\mathbb{P}- a.s.$ and satisfies
       \begin{equation}\label{L1}
    \lim\limits_{n\to\infty}c_n \in L^1 \quad\text{(resp}.\;\lim\limits_{n\to\infty}\chi_n\in L^1).
    \end{equation}
 \subsection{Preliminary results}
\subsubsection{Viscosity solutions }
\no In this part, we first give the definition of viscosity solution and then give some properties of the value function. \\
 Let $F$ be a  continuous function so that
 \begin{equation}\label{eq}
 F(x,v,D_{x}v,D^2_{xx}v)=0 ,
 \end{equation}
 where $D_x$ and $D_{xx}$ stand respectively  for the first derivative and the second derivative with respect to $x$. 
We assume moreover that $F$ is nonincreasing with respect to its last argument and $x$ is an element of $\mathcal{O}$ an open subset of $\mathbb{R}^+$.
\begin{axiom} Viscosity solution
\begin{itemize}
\item $v$ is said to be a viscosity supersolution, if for any $\overline{x}$ $\in$ $\mathcal{O}$, and any $C^2$ function $\phi$ in a neighborhood of $\overline{x}$ , $\overline{x}$ is a local minimum of $v-\phi$ and 
$$F(\overline{x},v(\overline{x}),D_{x}\phi(\overline{x}),D^2_{xx}\phi(\overline{x})) \geq 0 .$$
\end{itemize}
\begin{itemize}
\item $v$ is said to be a viscosity subsolution, if for any $\overline{x}$ $\in$ $\mathcal{O}$, and any $C^2$ function $\phi$ in a neighborhood of $\overline{x}$ , $\overline{x}$ is a local maximum of $v-\phi$ and 
$$F(\overline{x},v(\overline{x}),D_{x}\phi(\overline{x}),D^2_{xx}\phi(\overline{x})) \leq 0 .$$
\end{itemize}

\begin{itemize}
\item $v$ is said to be a viscosity solution if it is both a viscosity supersolution and a viscosity subsolution. 
\end{itemize}
Equivalently, Let's  denote by $J^{2,+}V(x)$ (resp.  $J^{2,-}V(x))$, the superjets (resp.  the subjets) of\, $V(x)$, the set of pairs $(q,X)\in \mathbb{R}^+\times \mathbb{S}$ such that:
\begin{equation*}
\begin{array}{ll}
\qquad V(y)\leq V(x)+\langle q,y-x\rangle+\frac{1}{2}\langle X(y-x),y-x\rangle+o(|y-x|^2),
 \\
(\text{resp .}  V(y)\geq V(x)+\langle q,y-x\rangle+\frac{1}{2}\langle X(y-x),y-x\rangle+o(|y-x|^2)). 
\end{array}
\end{equation*}

\begin{itemize}
    \item $v$ is a viscosity supersolution to (\ref{eq}) if
    $$  F(x,v,q,X) \geq 0  \quad \forall (q,X) \in J^{2,-}V(x)). $$
\end{itemize}

\begin{itemize}
    \item $v$ is a viscosity subsolution to (\ref{eq}) if
    $$  F(x,v,q,X) \leq 0  \quad \forall (q,X) \in J^{2,+}V(x).$$
    
\end{itemize}
\end{axiom}

\subsubsection{\textbf{Existence  of the value function}}
\no In a two player game framework, we have two value functions, namely the upper value function and the lower value function. 
\begin{axiom} Upper and lower  value function\\
Let $\mathcal{A}^i$ (resp. $\mathcal{B}^j$) be a set of all admissible strategies for player I within regime $i$, (resp. Player II within regime $j$). A nonanticipative strategy for player I is a mapping 
 $$\vartheta: \bigcup_{j \in \mathcal{D}} \mathcal{B}^j \rightarrow\mathcal{A}^i, $$ such that for any stopping time $\tau$ and any $b$ , $b'$ $\in$  $\mathcal{B}^j$, with $b=b'$ on $[0,\tau]$, it holds that $\vartheta(b)  =\vartheta(b')$ on $[0,\tau]. $ Similarly, a nonanticipative strategy for player \textit{II} is a mapping  $$\varrho: \bigcup_{i \in \mathcal{D}} \mathcal{A}^i \rightarrow\mathcal{B}^j , $$ such that for any stopping time $\tau$ and any $b$, $b'$ $\in$  $\mathcal{B}^j$, with $b=b'$ on $[0,\tau]$, it holds that $\varrho(b)  =\varrho(b')$ on $[0,\tau]. $\\ let's denote by  $\Gamma^i$ (resp $\Delta^j$) the set of all nonanticipative strategies $\alpha$ (resp. $\beta$) for player I (resp player II).
The upper value of the game then is given by :
\begin{equation}
\overline{V}_{ij}=\inf_{\varrho \in \Delta^j} \sup_{\varphi \in \mathcal{A}^i} J(x,\varphi, \varrho(\varphi)) ,  
\end{equation}
and the lower value of the game is given by :
\begin{equation}
\underline{V}_{ij}=\sup_{\varphi \in \Gamma^i} \inf_{\beta \in \mathcal{B}^j}  J(x,\varphi(\beta) , \beta).
\end{equation}
\end{axiom}
We now give some properties of the value functions.
\bl
\textit{ The process $X^{x}$ satisfies the following estimate:
There exists a constant $\rho$
such that}
\begin{equation}\label{estimat1}
E[|X^{x}_t|]\leq e^{\rho t}(1+|x|),\qquad \forall t \geq 0.
\end{equation}
\el 
\bpf
   We refer the reader to \textcolor{blue}{\cite{DM}} for more detail about this estimate. 
\epf
\bl
 The costs $c_{ij}$ and $\chi_{ij}$ satisfy the following properties
\begin{equation}\label{eqn2}
 - \sum_{n=1}^{N} e^{-r\tau_n}c(\xi_{n-1},\xi_n) \leq \max_{k \in \mathcal{D}} (-c(i,k)),
\end{equation}
\begin{equation}\label{eqn3}
\min_{k \in \mathcal{D}} (\chi(i,k)) \leq  \sum_{m=1}^{N} e^{-r\rho_m}\chi  (\eta_{m-1},\eta_m).
\end{equation}
\el

\bpf
First, Let's prove that:
$- \sum_{n=1}^{N} e^{-r\tau_n}c(\xi_{n-1},\xi_n) \leq \max_{k \in \mathcal{D}} (-c(i,k))$.
Using Induction, we can see that the result hold for $N=1$. Suppose it holds for some N-1, let us show that it also holds for N.
Le $A$ be defined as follows $A:=\{ \omega \in \Omega|C(\xi_{N-1},\xi_N) <0 \}$.
we have:
\begin{equation*}
\begin{array}{ll}
-\sum_{m=1}^{N}e^{-r\tau_m}c(\xi_{m-1},\xi_m)\\ \\ \qquad\qquad=\big(-\sum_{m=1}^{N-1}e^{-r\tau_m}Cc\xi_{m-1},\xi_m)-e^{-r\tau_{N}}Cc\xi_{N-1},\xi_{N}))\big)\mathbf{1}_{A^c}\\ \\ \qquad\qquad+\big(-\sum_{m=1}^{N-2}e^{-r\tau_m}C(\xi_{m-1},\xi_m)-(e^{-r\tau_{N-1}}C(\xi_{N-2},\xi_{N-1})+e^{-r\tau_{N}}Cc\xi_{N-1},\xi_{N}))\big)\mathbf{1}_{A}\\ \\
\qquad\qquad\leq \big(-\sum_{m=1}^{N-1}e^{-r\tau_m}Cc\xi_{m-1},\xi_m)\big)\mathbf{1}_{A^c}\\ \\ \qquad\qquad+\big(-\sum_{m=1}^{N-2}e^{-r\tau_m}c(\xi_{m-1},\xi_m)-e^{-r\tau_{N-1}}(c(\xi_{N-2},\xi_{N-1})+c(\xi_{N-1},\xi_{N}))\big)\mathbf{1}_{A}
\end{array}
\end{equation*}
\begin{equation*}
\begin{array}{ll}
\qquad\qquad\leq\big(-\sum_{m=1}^{N-1}e^{-r\tau_m}c(\xi_{m-1},\xi_m)\big)\mathbf{1}_{A^c}\\ \\
\qquad\qquad+\big( \sum_{l\in\mathcal{I}}\big [-\sum_{m=1}^{N-2}e^{-r\tau_m}c(\xi_{m-1},\xi_m)-e^{-r\tau_{N-1}}c(\xi_{N-2},l)\big]\mathbf{1}_{[\xi_{N}=l]}\big)\mathbf{1}_{A}\\ \\
\qquad\qquad\leq \big(-\sum_{m=1}^{N-1}e^{-r\tau_m}c(\xi_{m-1},\xi_m)\big)\mathbf{1}_{A^c}\\ \\
\qquad\qquad+\big(\sum_{l\in\mathcal{I}}\big[-\sum_{m=1}^{N-1}e^{-r\tau_m}c(\tilde{\xi}_{m-1},\tilde{\xi}_m)\big]\mathbf{1}_{[\xi_{N}=l]}\big)\mathbf{1}_{A}\leq \max_{k\in\mathcal{I}}(-c(i,k)) .
\end{array}
\end{equation*}
where  $\tilde{\xi}_m=\xi_m$ for $m=1,...,N-2,$ and $\tilde{\xi}_{N-1}=l$.\\
Now let's look at $ \min_{k \in \mathcal{D}} (\chi(i,k)) \leq  \sum_{m=1}^{N} e^{-r\rho_m}\chi  (\eta_{m-1},\eta_m)$. 
We have:\\ $\max_{k \in \mathcal{D}} (-c(i,k))=-\min_{k \in \mathcal{D}} (c(i,k))$, hence: $- \sum_{n=1}^{N} e^{-r\tau_n}C(\xi_{n-1},\xi_n) \leq -\min_{k \in \mathcal{D}} (c(i,k))$.
And finally multiplying both sides by -1 and replacing $c(i,j)$ by $\chi(i,j)$, we get
the desired result. \hfill $\Box$
\epf
\bl
Under assumptions (\textbf{H1}), (\textbf{H2}) and (\textbf{H3}) there exists some positive constant $\rho$ such that for $r \geq \rho$ , the lower and upper value functions satisfy a linear growth condition: for all $i$, $j$ $\in \mathcal{D}$, and $x \in \mathbb{R}^+$, there exist a constant C such that :
\begin{equation}
\label{lp}
|\overline{V}_{ij}(x)|, |\underline{V}_{ij}(x)| \leq C ( 1 +|x|).
\end{equation}
\el

\bpf
 Using (\ref{eq2}) and consider the particular strategy for player II:  $\tilde{\beta}:=(\tilde{\rho}_n,\tilde{\eta}_n)$  given by  $\tilde{\rho}_n=\infty, \tilde{\eta}_n=j$ for all $n\geq 1$ and by admissibility condition (\ref{L1}) for every $\delta\in\mathcal{A}^i$ we have
\begin{equation*}
 J(x,\delta,\tilde{\beta}(\delta) ) \leq \mathbb{E}\bigg[\integ{0}{\infty}e^{-rs}f^{a_sb_s}(X_{s}^{x,\delta,\tilde{\beta}(\delta) })ds+\max_{k\in\mathcal{I}}(-C(i,k))\bigg].
\end{equation*}
By the estimate (\ref{estimat1}) and the polynomial growth condition of $f^{ij}$  in (\textbf{H2}), there exists $C$ and $\rho$ such that
\begin{equation*}
 J(x,\delta,\tilde{\beta}(\delta) )\leq \frac{C}{r}+C(1+|x|)\integ{0}{\infty}e^{(\rho-r)s}ds+ \max_{k\in\mathcal{I}}(-C(i,k)).
\end{equation*}
Therefore from the arbitrariness of $\delta\in\mathcal{A}^i$ and  if $r>\rho$ we have
 \begin{equation}
 \overline{V}_{ij}(x)\leq C(1+|x|).
 \end{equation}
On the other hand, given that we have:
\begin{equation*}
\sum_{n=1}^{N}e^{-r\rho_n}\chi(\eta_{n-1},\eta_n)\geq \min_{l\in\mathcal{I}}(\chi(j,l)) \qquad a.s.
\end{equation*}
for all $N\geq 1.$ By considering
the particular strategy $\tilde{\delta}:=(\tilde{\tau}_m,\tilde{\xi}_m)$ given by $\tilde{\tau}_m=\infty,  
 \tilde{\xi}_m=i$ for all $m\geq 1$, and by admissibility condition (\ref{L1}) for every $\beta\in \Delta^j$ we have:
\begin{equation*}
J(x,\tilde{\delta},\beta(\tilde{\delta})) \geq \mathbb{E}\bigg[\integ{0}{\infty}e^{-rs}f^{a_sb_s}(X_{s}^{x,\tilde{\delta},\beta(\tilde{\delta})})ds+\min_{l\in\mathcal{I}}(\chi(j,l))\bigg].
\end{equation*}
As a consequence, there exists $C$ and $\rho$ such that:
\begin{equation*}
J(x,\tilde{\delta},\beta(\tilde{\delta})) \geq -\frac{C}{r}-C(1+|x|)\integ{0}{\infty}e^{(\rho-r)s}ds+\min_{l\in\mathcal{I}}(\chi(j,l)).
\end{equation*}
Therefore from the arbitrariness of $\beta\in\Delta^j$ and if $r>\rho$, we have:
$$\overline{V}_{ij}(x)\geq -C(1+|x|),$$
from which we deduce the  claim.  \hfill $\Box$
\epf 
The game is said to admit a value if we have $\overline{V}_{ij}=\underline{V}_{ij}$. This result has been proved in \textcolor{blue}{\cite{BM}} as well as the continuity of the value function.
\section{System of quasi-variational inequalities and switching regions }
\setcounter{equation}{0}
\renewcommand{\theequation}{3.\arabic{equation}}
We now state the Isaacs system of quasi-variational inequalities. Using the dynamic programming principle, it has been proved that  the problem of two player switching game leads to a systems of quasi-variational inequalities :
\begin{equation}
\label{qv}
\begin{array}{rr}
    max \bigg\{ min \bigg [ r { v}_{ij}(x)-\mathcal{L}_{ij} {v}_{ij}(x) - f_{ij}; v_{ij}(x)-M_{ij}[v](x)\bigg ]; v_{ij}(x)-N_{ij}[v](x) \bigg \}=0,
    \end{array}
\end{equation},
\begin{equation}
\begin{array}{rr}
    min \bigg\{ max \bigg [ r{ v}_{ij} (x) -\mathcal{L}_{ij} {v}_{ij}(x) - f_{ij}; v_{ij}(x)-M_{ij}[v](x)\bigg ]; v_{ij}(x)-N_{ij}[v](x)\bigg \}=0 ,
    \end{array}
\end{equation}
where  $$M_{ij}[v](x)= \max_{k\ne i}\{v_{kj}(x)-c_{ik}\}, \qquad N_{ij}[v](x)= \min_{l\ne j}\{v_{il}(x)+\chi_{jl}\},$$ 
and
$$\mathcal{L}_{ij}v_{ij}(x)=\frac{1}{2} \mathrm{Tr}[\sigma_{ij} (\sigma_{ij})^* \nabla_x^2 v_{ij}(x)]+ \big < b_{ij} ,\nabla_x v_{ij}(x)\big > .$$

The viscosity supersolution (resp subsolution) associated to the first quasi-inequality is that continuous function $v_{ij}(x)$ such that for any $i,j\in\mathcal{I}$, for any $x_0\in\mathbb{R}^+$ and any function $\phi_{ij}\in C^{2}(\mathbb{R}^+)$ such that $\phi_{ij}(x_0)=v_{ij}(x)(x_0)$ and $x_0$ is a local maximum (resp local minimum)
of $\phi_{ij}-v_{ij}(x)$, we have:
\begin{equation}
\begin{array}{lll}
max\Big\{ min\Big[r\phi_{ij}(x_0)-\mathcal{L}_{ij}\phi_{ij}(x_0)-f_{ij}(x_0),\\ \qquad v_{ij}(x)(x_0)-M_{ij}[v](x_0)\Big],v_{ij}(x)(x_0)-N_{ij}[v](x_0) \Big\}\geq 0,
\end{array}
\end{equation}
\begin{equation}
 \begin{array}{lll}
\bigg (\text{ resp } max\Big\{ min\Big[r\phi_{ij}(x_0)-\mathcal{L}_{ij}\phi_{ij}(x_0)-f_{ij}(x_0),\\ \qquad v_{ij}(x)(x_0)-M_{ij}[v](x_0)\Big],v_{ij}(x)(x_0)-N_{ij}[v](x_0) \Big\}\leq 0. \bigg )
\end{array}
\end{equation} 
Under the theorem of existence of the  value of the game, the game admit a solution if the upper and lower value function coincide.
\bl
    The upper and lower value function are solutions in the viscosity sense of both above quasi-variational inequality and both value functions coincide.
\el

\bpf
We refer the reader to \textcolor{blue}{\cite{BM}}  for the proof. \hfill $\Box$
\epf

\no We will now introduce the notion of switching regions and continuation regions. 
The switching regions correspond to the set of values of the diffusion process for which the optimal solution would be to change the regime for a given player.  Conversely, the continuation regions constitute the set of values of the process for which it would be optimal to stay in the same regime.

For $\xi - player$: \qquad \qquad  
\begin{equation*}
\begin{cases}
 S^{\xi}_{i}. &=\big \{x \in (0, \infty) : V_{i}.(x) = \max \limits_{k\in \mathcal{D}^{-i}} (V_{k}.(x) -c_{ik}) \big \},\\
 C^{\xi}_{i}. &=\big \{x \in (0, \infty) : V_{i}.(x)  > \max\limits_{k\in \mathcal{D}^{-i}} (V_{k}.(x) -c_{ik}) \big \}.
 \end{cases}
\end{equation*}
$ S^{\xi}_{i}. $ is a closed subset of $(0, \infty )$ and corresponds to the regime where it is optimal for player I ($\xi-player$) to switch to another regime, while its complement $ C^{\xi}_{i}.$  is an open subset of $(0, \infty )$.
\\
 For $\eta-player:$ \qquad \qquad \qquad
 \begin{equation*}
\begin{cases}
 S^{\eta}_{j}.  &=\big \{x \in (0, \infty) : V. _{j}(x)= \min\limits_{l\in \mathcal{D}^{-j}} (V. _{l}(x)+\chi_{jl}) \big \},\\
 C^{\eta}_{j}.  & =\big \{x \in (0, \infty) : V. _{j} (x)< \min \limits_{l\in \mathcal{D}^{-j}} (V. _{l} (x)+\chi_{jl}) \big \}.\\
\end{cases}
 \end{equation*}
\\
 $S^{\eta}_{j}. $ is a closed subset of $(0, \infty )$ and corresponds to the regime where it is optimal for player II ($\eta-player$) to switch to another regime, while its complement $C^{\eta}_{j}.$  is an open subset of $(0, \infty )$.
 \begin{rem}
The smooth-fit property highlighted in \textcolor{blue}{\cite{P}} is still valid in our context; under this property the value function is $C^1$ on $\partial S^{\eta}_{j}.$ and $\partial  S^{\xi}_{i}.$ and  $C^2$ on $C^{\eta}_{j}.$ and $C^{\xi}_{i}.$ .
 \end{rem}
 \bl
 \label{lem:lemma3}
In order for ${ v}_{ij}(x)$ to be solution of (\ref{qv}) one of the below conditions must be satisfied.
\begin{itemize}
    \item \textbf{Condition A1:}
    \begin{enumerate}[start=11,label={(\bfseries A\arabic*):}]
        \item  $r { v}_{ij}(x)-\mathcal{L}_{ij} {v}_{ij}(x) - f_{ij}(x)=0$;
        \item   $v_{ij}(x)-M_{ij}[v](x)\geq 0$;
        \item   $v_{ij}(x)-N_{ij}[v](x)\leq 0.$
    \end{enumerate} 
    \item \textbf{Condition A2:}
    \begin{enumerate}[start=21,label={(\bfseries A\arabic*):}]
        \item  $r { v}_{ij}(x)-\mathcal{L}_{ij} {v}_{ij}(x) - f_{ij}(x)\geq 0$;
        \item  $v_{ij}(x)-M_{ij}[v](x)= 0$;
        \item  $v_{ij}(x)-N_{ij}[v](x)\leq 0.$
    \end{enumerate}
    \item \textbf{Condition A3:}
   \begin{enumerate}[start=31,label={(\bfseries A\arabic*):}]
       \item  $v_{ij}(x)-N_{ij}[v](x)= 0;$
        \item  $min \big [ r v_{ij}(x)-\mathcal{L}_{ij} v_{ij}(x) - f_{ij}(x); v_{ij}(x)-M_{ij}[v](x)\big ]\leq 0 $.
    \end{enumerate} 
\end{itemize}
\el 
\bpf 
\begin{itemize}[leftmargin=*, itemsep=1.5 pt]
\item If condition A1 is verified:
  \begin{fleqn}[\parindent]
\begin{equation*}
   \textbf{(A11) + (A12)} \implies min \bigg [ r { v}_{ij}(x)-\mathcal{L}_{ij} {v}_{ij}(x) - f_{ij}(x); v_{ij}(x)-M_{ij}[v](x)\bigg ] =r { v}_{ij}(x)-\mathcal{L}_{ij} {v}_{ij}(x) - f_{ij}(x)=0; 
\end{equation*} 
\end{fleqn}
  \begin{fleqn}[\parindent]
  \begin{spreadlines}{2ex}
 \begin{equation*}
 \begin{split}
      \textbf{ \bigg [ (A11) + (A12) \bigg ]+ (A13) }  \implies & max \bigg \{ min \bigg [ r { v}_{ij}(x)-\mathcal{L}_{ij} {v}_{ij}(x) - f_{ij}(x); v_{ij}(x)-M_{ij}[v](x)\bigg ]; v_{ij}(x)-N_{ij}[v](x) \bigg \} \\
  &= r { v}_{ij}(x)-\mathcal{L}_{ij} {v}_{ij}(x) - f_{ij}(x)=0.
  \end{split}
 \end{equation*}
 \end{spreadlines}
 \end{fleqn}
\end{itemize}
\noindent Therefore, $v_{ij}(x) \quad i,j=1,2$ solves the system of quasi-variational inequalities (\ref{qv}).
\begin{itemize}[leftmargin=*, itemsep=1.5 pt]
 \item If condition A2 is verified:
   \begin{fleqn}[\parindent]
\begin{equation*}
\textbf{(A21) + (A22)} \implies  min \bigg [ r { v}_{ij}(x)-\mathcal{L}_{ij} {v}_{ij}(x) - f_{ij}(x); v_{ij}(x)-M_{ij}[v](x)\bigg ]=v_{ij}(x)-M_{ij}[v](x); 
\end{equation*}
\end{fleqn}
\begin{spreadlines}{2ex}
\begin{equation*}
\begin{split}
   \textbf{\bigg [(A21) + (A22) \bigg ] + (A23) } \implies & max \bigg\{ min \bigg [ r { v}_{ij}(x)-\mathcal{L}_{ij} {v}_{ij}(x) - f_{ij}(x); v_{ij}(x)-M_{ij}[v](x)\bigg ]; v_{ij}(x)-N_{ij}[v](x) \bigg \} \\ & = v_{ij}(x)-M_{ij}[v](x)=0. 
\end{split}
\end{equation*}
\end{spreadlines}
\noindent Therefore, $v_{ij}(x) \quad i,j=1,2$ solves the system of quasi-variational inequalities (\ref{qv}).
       \item If condition A3 is verified:
         \begin{fleqn}[\parindent]
         \begin{spreadlines}{2ex}
       \begin{equation*}
       \begin{split}
       \textbf{(A31) + (A32) } \implies & max \bigg\{ min \bigg [ r { v}_{ij}(x)-\mathcal{L}_{ij} {v}_{ij}(x) - f_{ij}(x); v_{ij}(x)-M_{ij}[v](x)\bigg ]; v_{ij}(x)-N_{ij}[v](x) \bigg \}\\
 & =  v_{ij}(x)-N_{ij}[v](x)=0.
       \end{split}
       \end{equation*}
       \end{spreadlines}
\end{fleqn}

\end{itemize}   
\noindent Therefore, $v_{ij}(x) \quad i,j=1,2$ solves the system of quasi-variational inequalities (\ref{qv}). \ \hfill $\Box$
\epf
\section{Explicit solution: Identical profit functions with different diffusion operators}
\setcounter{equation}{0}
\renewcommand{\theequation}{4.\arabic{equation}}
\no We are now about to give the explicit solution to the problem in the two regimes case ($d=2$). We have seen that the condition of existence of the value of the game implies that the upper value function and the lower value function coincide. So in order to solve our problem, we just solve the quasi-variational inequalities associated to the upper value function.   The 
solution of the problem are $v_{ij}(x)$ that satisfies the 4 quasi-variational inequalities below:
\begin{align}
 max \big \{ min\big [ r v_{11}(x)-\mathcal{L}_{11} v_{11}(x) - f_{11}(x); v_{11}(x)-(v_{21}(x)-c_{12}) \big ]; {v}_{11}(x)-(v_{12}(x) + \chi_{12}) \big \}=0,\label{eq1} \\
max \big \{ min \big [ r v_{12}(x)-\mathcal{L}_{12} v_{12}(x) - f_{12}(x); v_{12}(x)-(v_{22}(x)-c_{12}) \big ]; v_{12}(x)-(v_{11}(x) + \chi_{21}) \big \}=0, \label{eq2}\\
max \big \{ min \big [ r v_{21}(x)-\mathcal{L}_{21} v_{21}(x) - f_{21}(x); v_{21}(x)-(v_{11}(x)-c_{21}) \big ]; v_{21}(x)-(v_{22}(x) + \chi_{12}) \big \}=0, \label{eq3}\\
 max \big \{ min \big [ r v_{22}(x)-\mathcal{L}_{22} v_{22}(x) - f_{22}(x); v_{22}(x)-(v_{12}(x)-c_{21}) \big ]; v_{22}(x)-(v_{21}(x) + \chi_{21}) \big \} = 0. \label{eq4}
\end{align}
In the next subsections, we provide the  explicit solution in the case  of identical profit functions:
\begin{equation}
    f_{ij}(x)=x^{\gamma}, \qquad  0 < \gamma < 1,
\end{equation}with different diffusion operators.
\begin{rem}
    This function can sufficiently represent real-life profit function where profit grows at exponential rate with respect to a variable.
\end{rem}
 For each player, the decision to switch or not will be based on the final reward after switching compared to that of not switching.
Let us consider the following second order Ordinary Differential Equation (ODE):
\begin{equation}
\label{ODE}
    r v(x)-\mathcal{L}_{ij} v(x) -f_{ij}(x) =0.
\end{equation}
The system (\ref{qv}) boils down to this equation when both players are in their continuation region. The general solution without second member of this equation is given by:\\ 
$$ v(x)=A_{ij}x^{m^+_{ij}} +B_{ij} x^{m^-_{ij}},$$
where $A_{ij}$, $B_{ij}$ are some constants and, 
\begin{equation}
\label{mu}
m_{ij}^+ = -\frac{b_{ij}}{\sigma_{ij}^2} + \frac{1}{2} +  \sqrt{(-\frac{b_{ij}}{\sigma_{ij}^2} + \frac{1}{2})^2+ \frac{2r}{\sigma_{ij}^2}}  >1,
\\
\end{equation}
\begin{equation}
\label{md}
m_{ij}^-= -\frac{b_{ij}}{\sigma_{ij}^2} + \frac{1}{2}-  \sqrt{(-\frac{b_{ij}}{\sigma_{ij}^2} + \frac{1}{2})^2+ \frac{2r}{\sigma_{ij}^2}} <0.
\end{equation}
With regards to  (\ref{lp}), (\ref{mu}), (\ref{md}), we guess that if both players are in their continuation region in the neighborhood of $0$, then $B_{ij}=0$ and if they are in their continuation region in the neighborhood of $+ \infty$, then $A_{ij}=0.$\\
Let's also denote: 
$$\hat{V}_{ij}=E \bigg[ \int_0^{\infty} e^{-rt} f_{ij} (X^{x,ij}_t )dt \bigg ]. $$ Such
$\hat{V}_{ij}$ is a particular solution to ODE (\ref{ODE}).  It corresponds to the reward function associated with the non-switching strategy from the initial state $(x, i,j)$. A straightforward calculus shows that :
$$\hat{V}_{ij}=K_{ij}x^{\gamma}, \qquad\text{with } \qquad K_{ij}= \frac{  1}{ { r-b_{ij} \gamma +\frac{1}{2}\sigma_{ij}^2 \gamma (1-\gamma)}}< 0, \quad i,j=1,2. $$
We show that the structure of the switching regions actually depends only on the sign
of  $\text{K}_{i_1 j_1} - \text{K}_{i_2j_2}$, and of the sign of the switching costs $c_{ij}$ and $\chi_{ij}$. We will restrict ourselves in the case where the diffusion operators are equal to two by two ( $\text{K}_{i_1 j_1}=\text{K}_{i_2 j_2}$ ). Specific cases where  $\text{K}_{i_1i_1}=\text{K}_{j_1 j_1}$ $<$ $\text{K}_{i_2 i_2}=\text{K}_{j_2 j_2}$ or $\text{K}_{i_1i_1}=\text{K}_{j_1 j_1}$ $>$ $\text{K}_{i_2 i_2}=\text{K}_{j_2 j_2}$ where $i_1 \neq j_1$ and $i_2 \neq  j_2$ are not covered .
We assume moreover that $\text{K}_{i_1 j_1}=\text{K}_{i_2 j_2}$ is equivalent to $m_{i_1j_1}^+=m_{i_2j_2}^+$ and $m_{i_1j_1}^-=m_{i_2j_2}^-$.
We will be working with the following structure of switching costs:
\begin{table}[H]
\centering
\begin{tabular}{|l|l|l|l|l|l|l|}
   \hline
& $c_{12}$ & $c_{21}$ & $\chi_{12}$ &   $\chi_{21}$ & $c_{12}$ +  $c_{21}$ & $\chi_{12}$ + $\chi_{21}$  \\
\hline
\hline
   \textbf{Condition B1} & \quad > 0 & \quad > 0 & \quad > 0 & \quad > 0 & \quad > 0 & \quad > 0\\
   \hline
  \textbf{Condition B2}   & \quad < 0 & \quad > 0 & \quad > 0 & \quad > 0 & \quad > 0 & \quad > 0 \\ 
  \hline
  \textbf{Condition B3}   & \quad > 0 & \quad > 0 & \quad < 0 & \quad > 0 & \quad > 0 & \quad > 0 \\
  \hline
   \textbf{Condition B4}   & \quad < 0 & \quad > 0 & \quad < 0 & \quad > 0 & \quad > 0 & \quad > 0\\
   \hline
 \end{tabular}
 \end{table} 
 \begin{rem}
Although these conditions and assumptions above might seem restrictive, they lead to 20 systems of quasi-variational inequalities to be solved, that is, 80 quasi-variational inequalities. This is quite enough as first approach on explicit solution of a two-player switching game problem. We can also use  symmetrical reasoning to get solutions for other structure of switching costs.
 \end{rem}
\bethe
\normalfont
\label{th:theo4.1}
Case where $K_{11} =   K_{12}  =  K_{21}  =  K_{22}$.
\begin{enumerate}[leftmargin=*, itemsep=1.5 pt]
    \item  If  Condition B1 is verified, then  let $v_{ij}(x)$ be defined as follows:
    \begin{tcolorbox}[colback=yellow!10!white, colframe=red!50!black]
    \begin{fleqn}[\parindent]
 \begin{equation*}
     \eqskip v_{ij}(x)=\hat{V}_{ij}(x),\quad x \in (0,+ \infty) \quad  i,j=1,2 .
 \end{equation*} 
    \end{fleqn}
    \end{tcolorbox}
    Such $v_{ij}(x)$ are the solutions associated to the system of quasi-variational inequalities with:
    
    \begin{center}
    $S^{\xi}_{ij}= \emptyset, \quad C^{\xi}_{ij}= (0, +\infty), \quad S^{\eta}_{ij}= \emptyset, \quad C^{\eta}_{ij}=(0, +\infty), \qquad i,j=1,2.$
   \end{center}

    It is never optimal for both players to switch when their respective switching costs are positive .
   \item If Condition B2 is verified, then let $v_{ij}(x)$ be defined as follows: 

\begin{tcolorbox}[colback=yellow!10!white, colframe=red!50!black]
\begin{minipage}{0.48\textwidth}
\begin{fleqn}[\parindent]
\begin{spreadlines}{1ex}
     \begin{align}
     v_{11}(x)&=\hat{V}_{11}(x)-c_{12},\\
     v_{12}(x)&=\hat{V}_{12}(x)-c_{12},
     \end{align}
     \end{spreadlines}
     \end{fleqn}
\end{minipage}
\hfill\vline\hfill
\begin{minipage}{0.5\textwidth}
\begin{fleqn}[\parindent]
\begin{spreadlines}{1ex}
   \begin{align}
   \; v_{21}(x)&=\hat{V}_{21}(x), \\
   \; v_{22}(x)&=\hat{V}_{22}(x).
     \end{align}
     \end{spreadlines}
     \end{fleqn}
\end{minipage}
\end{tcolorbox}

Such $v_{ij}(x)$ are the solutions associated to the system of quasi-variational inequalities with: 

\begin{center}
     $S^{\xi}_{11}= S^{\xi}_{12}=(0, +\infty), \quad
 S^{\eta}_{11}= S^{\eta}_{12}=S^{\xi}_{21}=S^{\eta}_{21}=S^{\xi}_{22}= S^{\eta}_{22}=\emptyset.$
\end{center}

When cost of switching from regime $1$ to regime $2$ is negative for player I (condition B2), it is always optimal for him to switch from regime $1$ to regime $2$. Player II on his side has to never switch.
 \item If  Condition B3 is verified, then let $v_{ij}(x)$ be defined as follows: 

\begin{tcolorbox}[colback=yellow!10!white, colframe=red!50!black]
\begin{minipage}{0.48\textwidth}
\begin{fleqn}[\parindent]
\begin{spreadlines}{1ex}
     \begin{align}
v_{11}(x)&=\hat{V}_{11}(x)+\chi_{12}, \\
v_{12}(x)&=\hat{V}_{12}(x), 
     \end{align}
     \end{spreadlines}
     \end{fleqn}
\end{minipage}
\hfill\vline\hfill
\begin{minipage}{0.5\textwidth}
\begin{fleqn}[\parindent]
\begin{spreadlines}{1ex}
   \begin{align}
   \; v_{21}(x)&=\hat{V}_{21}(x)+ \chi_{12}, \\
   \;  v_{22}(x)&=\hat{V}_{22}(x).
     \end{align}
     \end{spreadlines}
     \end{fleqn}
\end{minipage}
\end{tcolorbox}

    Such $v_{ij}(x)$ are the solutions associated to the system of quasi-variational inequalities with:

\begin{center}
    $ S^{\xi}_{11}= S^{\xi}_{12}=S^{\eta}_{12}=S^{\xi}_{21}=S^{\xi}_{22}=S^{\eta}_{22}= \emptyset, \quad  S^{\eta}_{11}=S^{\eta}_{21}=(0, +\infty).$
\end{center}

When cost of switching from regime $1$ to regime $2$ is negative for player II (condition B3), it is always optimal for him to switch from regime $1$ to regime $2$. Player I on his side has to never switch.
   \item If Condition B4 is verified, then  let $v_{ij}(x)$ be defined as follows: 

\begin{tcolorbox}[colback=yellow!10!white, colframe=red!50!black]
\begin{minipage}{0.48\textwidth}
\begin{fleqn}[\parindent]
\begin{spreadlines}{1ex}
     \begin{align}
      v_{11}(x)&=\hat{V}_{11}(x)-c_{12} +\chi_{12}, \\
      v_{12}(x)& =\hat{V}_{12}(x)- c_{12},
     \end{align}
     \end{spreadlines}
     \end{fleqn}
\end{minipage}
\hfill\vline\hfill
\begin{minipage}{0.5\textwidth}
\begin{fleqn}[\parindent]
\begin{spreadlines}{1ex}
   \begin{align}
  \;  v_{21}(x)&=\hat{V}_{21}(x)+ \chi_{12},\\
   \;  v_{22}(x)&=\hat{V}_{22}(x).
     \end{align}
     \end{spreadlines}
     \end{fleqn}
\end{minipage}
\end{tcolorbox}

Such $v_{ij}(x)$ are the solutions associated to the system of quasi-variational inequalities with: 

  \begin{center}
$S^{\xi}_{11}=S^{\eta}_{11}= S^{\xi}_{12}=S^{\eta}_{21}=(0, +\infty), \quad   S^{\eta}_{12}= S^{\xi}_{21}= S^{\xi}_{22}= S^{\eta}_{22}=\emptyset.$
   \end{center}
   
   When cost of switching from regime $1$ to regime $2$ is negative for both players (condition B4),  it is always optimal for them to switch from regime $1$ to regime $2$.
\end{enumerate}    
\eethe

\bpf
\normalfont
\begin{enumerate}[leftmargin=*]
    \item  Consider the case of condition (B1) and let $v_{ij}(x)$ be defined as in Theorem \ref{th:theo4.1} -1:\\  From above, we know that $\hat{V}_{ij}(x)$ is a particular solution to (\ref{ODE}), 
    hence, $v_{ij}(x)$ satisfies (A11). 
     On the other hand, we have:
     $v_{ij}(x)- \big (v_{kj}(x)-c_{ik} \big )=\hat{V}_{ij}(x)-\big (\hat{V}_{kj}(x)-c_{ik} \big )=c_{ik} >0$. 
     Therefore $v_{ij}(x)$ satisfies (A12) as well.
     With regards to  $(A13)$, we have:
     $v_{ij}(x)-\big (v_{ik}(x) + \chi_{jk} \big )= \hat{V}_{ij}(x)-\big (\hat{V}_{kj}(x)-\chi_{jk} \big )=- \chi_{jk}<0.$
     By Lemma \ref{lem:lemma3}, we deduce that $v_{ij}(x)$ solves the system \ref{qv}. \hfill  $\Box$
    \item  Consider the case of condition (B2): Let $v_{ij}(x)$ be defined as in Theorem \ref{th:theo4.1} -2:
    \begin{itemize}[label=$\star$, leftmargin=*]   \setlength\itemsep{1em}
        \item  We have by definition of $v_{ij}(x)$:
        $rv_{11}(x)-\mathcal{L}_{11} v_{11}(x) - f_{11}(x)=-rc_{12}>0,$ hence $v_{ij}(x)$ satisfies (A21). On the other hand, we have: 
        $v_{11}(x)-\big (v_{21}(x)-c_{12} \big )=\hat{V}_{11}(x)-c_{12} - \big (\hat{V}_{21}(x)-c_{12} \big )=0,$
        therefore $v_{ij}(x)$ satisfies (A22) as well. With regards to (A23), we have: $v_{11}(x)- \big (v_{12}(x) + \chi_{12} \big )=\hat{V}_{11}(x)-c_{12}-\big (\hat{V}_{12}(x)-c_{12} + \chi_{12} \big )= -\chi_{12}<0.$ By Lemma \ref{lem:lemma3}, we deduce that $v_{ij}(x)$ solves (\ref{eq1}).
        \item We have by definition of $v_{ij}(x)$ : $rv_{12}(x)-\mathcal{L}_{12} v_{12}(x) - f_{12}(x)=-rc_{12}>0,$ hence $v_{ij}(x)$ satisfies (A21). On the other hand,  $v_{12}(x)-\big (v_{22}(x)-c_{12} \big )=\hat{V}_{12}(x)-c_{12} - \big (\hat{V}_{22}(x)-c_{12} \big )=0,$ thus, $v_{ij}(x)$ satisfies (A22) as well. With regards to (A23), we have $v_{12}(x)-\big (v_{22}(x) + \chi_{12} \big )=\hat{V}_{12}(x)-c_{12}-\big (\hat{V}_{22}(x)-c_{12} + \chi_{12} \big )= -\chi_{12}<0.$ By Lemma \ref{lem:lemma3}, we deduce that $v_{ij}(x)$ solves (\ref{eq2}).
        \item We have by definition of $v_{ij}(x)$:
        $rv_{21}(x)-\mathcal{L}_{21} v_{21}(x) - f_{21}(x)= 0,$ therefore, $v_{ij}(x)$ satisfies (A11). On the other hand, we have: $v_{21}(x)-\big (v_{11}(x)-c_{21} \big )=\hat{V}_{21}(x)- \big (\hat{V}_{11}(x) -c_{12}-c_{21} \big )=(c_{12}+c_{21})>0,$ thus, $v_{ij}(x)$ satisfies  (A12). With regards to (A13) we have: $v_{21}(x)-\big (v_{22}(x) + \chi_{21} \big ) =\hat{V}_{21}(x)-\big (\hat{V}_{22}(x) + \chi_{21} \big ) = - \chi_{21} <0.$ By Lemma \ref{lem:lemma3}, we deduce that $v_{ij}(x)$ solves(\ref{eq3}).
        \item We have by definition of $v_{ij}(x)$ : $r v_{22}(x)-\mathcal{L}_{22} v_{22}(x) - f_{22}(x)=0,$ thus, $v_{ij}(x)$ satisfies (A11). On the other hand, we have: $v_{22}(x)-\big (v_{12}(x)-c_{21} \big )=\hat{V}_{22}(x)-\big (\hat{V}_{12}(x) -c_{12}-c_{21} \big )=(c_{12}+c_{21})>0,$ hence $v_{ij}(x)$ satisfies (A12). With regards to (A13), we have: $v_{22}(x)-\big (v_{21}(x) + \chi_{21} \big ) =\hat{V}_{22}(x)-\big (\hat{V}_{21}(x) + \chi_{21} \big )=- \chi_{21} <0.$  By Lemma \ref{lem:lemma3}, we deduce that $v_{ij}(x)$ solves(\ref{eq4}).  \hfill $\Box$
    \end{itemize}
    \item Consider the case condition (B3): Let $v_{ij}(x)$ be defined as in Theorem \ref{th:theo4.1} -3:
    
    \begin{itemize}[label=$\star$, leftmargin=*]   \setlength\itemsep{1em}
        \item  We have by definition of $v_{ij}(x)$:
        $rv_{11}(x)-\mathcal{L}_{11} v_{11}(x) - f_{11}(x)=r\chi_{12}<0,$ thus $v_{ij}(x)$ satisfies (A32). On the other hand, we have:  $v_{11}(x)-\big (v_{12}(x) + \chi_{12} \big )=\hat{V}_{11}(x)+\chi_{12}- \big (\hat{V}_{12}(x)+\chi_{12} \big )= 0,$ hence $v_{ij}(x)$ satisfies (A31). By Lemma \ref{lem:lemma3}, we deduce that $v_{ij}(x)$ solves (\ref{eq1}).
        \item We have  by definition of $v_{ij}(x)$:
        $rv_{12}(x)-\mathcal{L}_{12} v_{12}(x) - f_{12}(x)=r\chi_{12}<0,$ thus, $v_{ij}(x)$ satisfies (A32). On the other hand, we have: $v_{12}(x)-\big (v_{22}(x) + \chi_{12} \big )=\hat{V}_{12}(x)+ \chi_{12}- \big (\hat{V}_{22}(x)+ \chi_{12} \big )= 0,$ therefore, $v_{ij}(x)$ satisfies (A31). By Lemma \ref{lem:lemma3}, we deduce that $v_{ij}(x)$ solves (\ref{eq2}).
        \item  We have by definition of $v_{ij}(x)$:  $rv_{21}(x)-\mathcal{L}_{21} v_{21}(x) - f_{21}(x)=0,$ thus, $v_{ij}(x)$ satisfies (A11). On the other hand, $v_{21}(x)-\big (v_{11}(x)-c_{21} \big )=\hat{V}_{21}(x)-\big (\hat{V}_{11}(x) +\chi_{12} -c_{21} \big )  = c_{21} - \chi_{12}>0,$ hence $v_{ij}(x)$ satisfies (A12). With regards to (A13) we have: $v_{21}(x)- \big (v_{22}(x) + \chi_{21} \big ) =\hat{V}_{21}(x)- \big (\hat{V}_{22}(x) + \chi_{21} \big )=-\chi_{21}  <0.$  By Lemma \ref{lem:lemma3}, we deduce that $v_{ij}(x)$ solves (\ref{eq3}).
        \item We have by definition of $v_{ij}(x)$: $r v_{22}(x)-\mathcal{L}_{22} v_{22}(x) - f_{22}(x)=0,$ hence $v_{ij}(x)$ satisfies (A11). On the other hand, we have: $v_{22}(x)-\big (v_{12}(x)-c_{21} \big )=\hat{V}_{22}(x)- \big (\hat{V}_{12}(x) +\chi_{12}-c_{21} \big )= c_{21} - \chi_{12}>0,$ hence $v_{ij}(x)$ satisfies (A12). With regards to (A13), we have: $v_{22}(x)- \big (v_{21}(x) + \chi_{21} \big ) =\hat{V}_{22}(x)- \big (\hat{V}_{21}(x) + \chi_{21} \big )=-\chi_{21}<0.$ By Lemma \ref{lem:lemma3}, we deduce that $v_{ij}(x)$ solves (\ref{eq4}). \hfill $\Box$
    \end{itemize}
    \item Consider the case of condition (B4):  Let $v_{ij}(x)$ be defined as in Theorem \ref{th:theo4.1} -4:
\begin{itemize}[label=$\star$, leftmargin=*]   \setlength\itemsep{1em}
        \item  We have by definition of $v_{ij}(x)$: $v_{11}(x)-\big (v_{21}(x)-c_{12} \big )=\hat{V}_{11}(x)-c_{12} + \chi_{12}- \big (\hat{V}_{21}(x)+ \chi_{12}-c_{12} \big )=0,$ thus, $v_{ij}(x)$ satisfies (A32). On the other hand, we have: $v_{11}(x)-\big (v_{12}(x) + \chi_{12} \big )=\hat{V}_{11}(x)-c_{12}+ \chi_{12}- \big (\hat{V}_{12}(x)-c_{12} + \chi_{12} \big )= 0,$ thus, $v_{ij}(x)$ satisfies (A31). By Lemma \ref{lem:lemma3}, we deduce that $v_{ij}(x)$ solves (\ref{eq1}).
        \item We have by definition of $v_{ij}(x)$:  $r v_{12}(x)-\mathcal{L}_{12} v_{12}(x) - f_{12}(x)=-rc_{12}>0,$ thus, $v_{ij}(x)$ satisfies (A21). On the other hand, we have:  $v_{12}(x)-\big (v_{22}(x)-c_{12} \big )=\hat{V}_{12}(x)-c_{12} -\big (\hat{V}_{22}(x)-c_{12} \big )=0,$ therefore, $v_{ij}(x)$ satisfies (A22). With regards to (A23), we have: $r v_{12}(x)-\big (v_{11}(x) + \chi_{21} \big )=\hat{V}_{12}(x)-c_{12}-\big (\hat{V}_{11}(x)-c_{12} + \chi_{21} \big )= - (\chi_{21} + \chi_{12})<0.$ By Lemma \ref{lem:lemma3}, we deduce that $v_{ij}(x)$ solves (\ref{eq2}).
        \item We have by definition of $v_{ij}(x)$: $r v_{21}(x)-\mathcal{L}_{21} v_{21}(x) - f_{21}(x)=r\chi_{12} < 0,$ thus, $v_{ij}(x)$ satisfies (A32). On the other hand, we have: $v_{21}(x)- \big (v_{22}(x) + \chi_{21} \big ) =\hat{V}_{21}(x)+ \chi_{12}- \big (\hat{V}_{22}(x) + \chi_{12} \big )=0.$ By Lemma \ref{lem:lemma3}, we deduce that $v_{ij}(x)$ solves (\ref{eq3}).
        \item We have by definition of $v_{ij}(x)$:  $r v_{22}(x)-\mathcal{L}_{22} v_{22}(x) - f_{22}(x)=0,$ thus,  $v_{ij}(x)$ satisfies (A11). On the other hand, we have: $v_{22}(x)-\big (v_{12}(x)-c_{21} \big )=\hat{V}_{22}(x)- \big (\hat{V}_{12}(x) -c_{12}-c_{21} \big )=(c_{12}+c_{21})>0,$ hence $v_{ij}(x)$ satisfies  (A12). With regards to (A13), we have: $v_{22}(x)- \big (v_{21}(x) + \chi_{21} \big ) =\hat{V}_{22}(x)- \big (\hat{V}_{21}(x) + \chi_{21}+\chi_{12} \big )=-(\chi_{21}+\chi_{12})<0.$ By Lemma \ref{lem:lemma3}, we deduce that $v_{ij}(x)$ solves (\ref{eq4}). \hfill $\Box$
    \end{itemize}
\end{enumerate}
    \epf
  \bethe
  \normalfont
  \label{th:theo4.2}
    Case where $K_{11} = K_{12}  <  K_{21} = K_{22}.$
    \begin{enumerate}[leftmargin=*]
        \item If Condition B1 is verified, then let $v_{ij}(x)$ be defined as follows
\begin{tcolorbox}[colback=yellow!10!white, colframe=red!50!black]
\begin{minipage}{0.48\textwidth}
\begin{fleqn}[\parindent]
\begin{spreadlines}{1ex}
     \begin{align}
     v_{11}(x)&= 
     \begin{cases}
      \hat{V}_{11}(x) +  Ax^{m^+_{11}}, \;  &x < x^*\\
      \hat{V}_{21}(x)  -  c_{12}, \; & x\geq x^*
     \end{cases},
     \\
     \vspace{0.3 cm}
     v_{21}(x) & =  \hat{V}_{21}(x),   \qquad \qquad \qquad x >0,
     \end{align}
     \end{spreadlines}
     \end{fleqn}
\end{minipage}
\hfill\vline\hfill
\begin{minipage}{0.5\textwidth}
\begin{fleqn}[\parindent]
\begin{spreadlines}{1ex}
   \begin{align}
   \;  v_{12}(x) &= \begin{cases}
     \hat{V}_{12}(x) +  Ax^{m^+_{12}}, \; & x < x^*\\
     \hat{V}_{22}(x)  -  c_{12}, \; & x\geq x^*
    \end{cases}, 
    \\
    \vspace{ 0.3 cm}
  \;  v_{22}(x) & =  \hat{V}_{22}(x), \qquad \qquad \qquad   x > 0,
     \end{align}
     \end{spreadlines}
     \end{fleqn}
\end{minipage}
\end{tcolorbox}

where $x^*=\inf S^{\xi}_{12} \in (0, +\infty)$. Using the fact that $v$ is $\mathcal{C}^1$ on $\partial S^{\xi}_{12}$ we get :
$$(K_{21}-K_{11})(x^*)^{\gamma}=\frac{m_{11}^{+} c_{12}}{m_{11}^{+}- \gamma}, \quad
A=(K_{21}-K_{11})\frac{\gamma}{m_{11}^{+}}(x^*)^{\gamma -m_{11}^{+}}.$$

Such $v_{ij}(x)$ are the solutions associated to the system of quasi-variational inequalities with:
\begin{center}
$S^{\xi}_{11}=S^{\xi}_{12}=[ x^*, +\infty), \quad
 S^{\eta}_{11}= S^{\eta}_{12}=S^{\xi}_{21}= S^{\eta}_{21}=S^{\xi}_{22}=S^{\eta}_{22}=\emptyset. $
 \end{center}
When the switching costs are positive (condition B1), player I has to stay in regime 1 for $x$ less than the threshold value $x^*$ and switch beyond that value. Player II on his side has to never switch.
\item If Condition B2 is verified then let $v_{ij}(x)$ be defined as follows:

\begin{tcolorbox}[colback=yellow!10!white, colframe=red!50!black]
\begin{minipage}{0.48\textwidth}
\begin{fleqn}[\parindent]
\begin{spreadlines}{1ex}
     \begin{align}
     v_{11}(x)  =  \hat{V}_{21}(x)- c_{12}, \\
     v_{12}(x) =  \hat{V}_{22}(x)- c_{12}, 
     \end{align}
     \end{spreadlines}
     \end{fleqn}
\end{minipage}
\hfill\vline\hfill
\begin{minipage}{0.5\textwidth}
\begin{fleqn}[\parindent]
\begin{spreadlines}{1ex}
   \begin{align}
   \;  v_{21}(x)  =  \hat{V}_{21}(x),\\ 
    \vspace{ 0.3 cm}
  \;  v_{22}(x) =  \hat{V}_{22}(x).
     \end{align}
     \end{spreadlines}
     \end{fleqn}
\end{minipage}
\end{tcolorbox}
Such $v_{ij}(x)$ are the solutions associated to the system of quasi-variational inequalities with: 
\begin{center}
$S^{\xi}_{11}=S^{\xi}_{12}=(0, +\infty), \quad
S^{\eta}_{11}=S^{\eta}_{12}=S^{\xi}_{21}= S^{\eta}_{21}= S^{\xi}_{22}= S^{\eta}_{22}=\emptyset.$ \end{center}
When cost of switching from regime $1$ to regime $2$ is negative for player I (condition B2), it is always optimal for him to switch from regime $1$ to regime $2$.
 
\item If Condition B3 is verified, then let $v_{ij}(x)$ be defined as follows:
\begin{adjustwidth}{-1 cm}{-1 cm}
\begin{tcolorbox}[colback=yellow!10!white, colframe=red!50!black]
\begin{minipage}{0.54\textwidth}
\begin{fleqn}[\parindent]
\begin{spreadlines}{1ex}
     \begin{align}
     v_{11}(x) & = 
  \begin{cases}
\hat{V}_{11}(x) +  Ax^{m^+_{11}} + \chi_{12}, & x < x^*\\
\hat{V}_{21}(x)  -  c_{12}+ \chi_{12,} & x\geq x^*
\end{cases}, \\
\; v_{21}(x)&  =  \hat{V}_{21}(x)+ \chi_{12},  \qquad \qquad \qquad x >0, 
    \end{align}
     \end{spreadlines}
     \end{fleqn}
\end{minipage}
\hfill\vline\hfill
\begin{minipage}{0.44\textwidth}
\begin{fleqn}[\parindent]
\begin{spreadlines}{1ex}
   \begin{align}
\; v_{12}(x)& =\begin{cases}
\hat{V}_{12}(x) +  Ax^{m^+_{12}},&   x < x^*\\
\hat{V}_{22}(x)  -  c_{12}, &  x\geq x^*
\end{cases},\\
\; v_{22}(x)& =  \hat{V}_{22}(x),\qquad \qquad \qquad  x > 0,
     \end{align}
     \end{spreadlines}
     \end{fleqn}
\end{minipage}
\end{tcolorbox}
\end{adjustwidth}
 where $x^*=\inf S^{\xi}_{12} \in (0, +\infty)$. Using the fact that $v$ is $\mathcal{C}^1$ on $\partial S^{\xi}_{12}$ we get :
$$(K_{21}-K_{11})(x^*)^{\gamma}=\frac{m_{11}^{+} c_{12}}{m_{11}^{+}- \gamma}, \quad
A=(K_{21}-K_{11})\frac{\gamma}{m_{11}^{+}}(x^*)^{\gamma -m_{11}^{+}}.$$
Such $v_{ij}(x)$ are the solutions associated to the system of quasi-variational inequalities with:
$$S^{\xi}_{11}=S^{\xi}_{12}=[x^*,\infty),\quad
S^{\eta}_{11}=S^{\eta}_{21}=(0,+\infty),  \quad S^{\xi}_{22}=S^{\eta}_{22}=S^{\eta}_{12}=
    S^{\xi}_{21}=\emptyset.$$
When the switching cost is negative from regime $1$ to regime $2$ for player II (condition B3), player I has to stay in regime 1 for $x$ less than the threshold value $x^*$ and switch beyond that value, while player II has to switch from regime $1$ to regime $2$ no matter the value of the state variable $x$.
    \item If Condition B4 is verified, then let $v_{ij}(x)$ be defined as follows:

\begin{tcolorbox}[colback=yellow!10!white, colframe=red!50!black]
\begin{minipage}{0.48\textwidth}
\begin{fleqn}[\parindent]
\begin{spreadlines}{1ex}
     \begin{align}
    \; v_{11}(x) & =  \hat{V}_{21}(x)- c_{12} + \chi_{12}, \\
   \;  v_{12}(x) & =  \hat{V}_{22}(x)- c_{12}, 
     \end{align}
     \end{spreadlines}
     \end{fleqn}
\end{minipage}
\hfill\vline\hfill
\begin{minipage}{0.5\textwidth}
\begin{fleqn}[\parindent]
\begin{spreadlines}{1ex}
   \begin{align}
  \; v_{21}(x) & =  \hat{V}_{22}(x) + \chi_{12}, \\
  \; v_{22}(x) &=  \hat{V}_{22}(x).
     \end{align}
     \end{spreadlines}
     \end{fleqn}
\end{minipage}
\end{tcolorbox}

Such $v_{ij}(x)$ are the solutions associated to the system of quasi variational inequalities with:
$$S^{\xi}_{11}=S^{\eta}_{11}= S^{\xi}_{12}= S^{\eta}_{21}=(0,+\infty), \quad S^{\eta}_{12}=S^{\xi}_{21}= S^{\xi}_{22}=S^{\eta}_{22}=\emptyset.$$
When the switching cost from regime $1$ to regime $2$ is negative for both players (condition B4), it is always optimal for them to switch from regime $1$ to regime $2$. 
    \end{enumerate}
\eethe

\bpf
\normalfont
\begin{enumerate}[leftmargin=*]
    \item Consider the case of condition (B1) and let $v_{ij}(x)$ be defined as in Theorem \ref{th:theo4.2} -1;\\ let us first check solutions for $x <  x^{*} $:
\begin{itemize}[label=$\star$, leftmargin=*]   \setlength\itemsep{1em}
    \item We have by definition of  $v_{ij}(x)$.
     $r v_{11}(x)-\mathcal{L}_{11} v_{11}(x) - f_{11}(x)=0,$
    thus, $v_{ij}(x)$ satisfies (A11). On the other hand, we have: $v_{11}(x)-\big (v_{21}(x)- c_{12} \big )=\hat{V}_{11}(x)-\hat{V}_{21}(x)+ Ax^{m^+_{11}}+ c_{12}$. Let $g$ be defined as follows: $g(x)=\hat{V}_{11}(x)-\hat{V}_{21}(x)+ Ax^{m^+_{11}}+c_{12}= (K_{11} -K_{21}) x ^{\gamma}+ Ax^{m^+_{11}}+c_{12}.$
A straightforward calculus shows that $g'(x)$ is negative, this implies that g(x) is a decreasing function; moreover, $g(x^*)=0 $, as such $g(x)$ is positive for $x < x^*$. Therefore, $v_{ij}(x)$ satisfies (A12). With regards to (A13), we have: $v_{11}(x)-\big (v_{12}(x)+ \chi_{12} \big )=-\chi_{12}<0.$ By Lemma \ref{lem:lemma3}, we deduce that $v_{ij}(x)$ solves (\ref{eq1}).
 \item We have by definition of  $v_{ij}(x)$: 
  $r v_{21}(x)-\mathcal{L}_{21} v_{21}(x) - f_{21}(x)=0,$
 thus, $v_{ij}(x)$ satisfies (A11). On the other hand, we have: $v_{21}(x)- \big (v_{11}(x) - c_{21} \big )= -g(x)+ c_{21}+c_{12}$ where $g$ is the function from above. This implies that $v_{ij}(x)$ satisfies (A12). With regards to (A13) we have:  $v_{21}(x)-\big (v_{22}(x) +\chi_{21} \big )=- \chi_{21} <0.$ By Lemma \ref{lem:lemma3}, we deduce that $v_{ij}(x)$ solves (\ref{eq3}).
\end{itemize}
We can notice that $v_{11}(x)=v_{12}(x)$ and $v_{21}(x)=v_{22}(x)$, hence the proof for  (\ref{eq1}) applies to  (\ref{eq2}) and the proof for  (\ref{eq3}) applies to  (\ref{eq4}). 

Let us now check solutions for $x \geq  x^{*} $:
\begin{itemize}[label=$\star$, leftmargin=*]   \setlength\itemsep{1em}
    \item  A straightforward calculus shows that by definition of $v_{ij}(x)$ we have:
    $r v_{11}(x)-\mathcal{L}_{11} v_{11}(x) - f_{11}(x)= \frac{ K_{21} -K_{11}}{K_{11}} x^{\gamma} - r c_{12}.$
The right hand side of the above equation is positive, hence $v_{ij}(x)$ satisfies (A21). On the other hand, we have : $v_{11}(x)-\big (v_{21}(x) + c_{12} \big )=\hat{V}_{21}(x) - c_{12} - \big ( \hat{V}_{21}(x)  - c_{12} \big )=0,$ thus, $v_{ij}(x)$ satisfies (A22). With regards to (A23), we have:  $v_{11}(x)-\big (v_{12}(x) + \chi_{12} \big )=0$. By Lemma \ref{lem:lemma3}, we deduce that $v_{ij}(x)$ solves (\ref{eq1}).
\item   We have by definition of  $v_{ij}(x)$: 
  $r v_{21}(x)-\mathcal{L}_{21} v_{21}(x) - f_{21}(x)=0,$
 thus, $v_{ij}(x)$ satisfies (A11). On the other hand, we have: $v_{21}(x)-\big (v_{11}(x) - c_{21} \big )=  c_{21}+c_{12}>0,$  hence $v_{ij}(x)$ satisfies (A12). With regards to (A13) we have: $v_{21}(x)-\big (v_{22}(x) +\chi_{21} \big )=- \chi_{21} <0$. By Lemma \ref{lem:lemma3}, we deduce that $v_{ij}(x)$ solves (\ref{eq3}).
\end{itemize}
Again, we see that $v_{11}(x)=v_{12}(x)$ and $v_{21}(x)=v_{22}(x)$, hence the proof for  (\ref{eq1}) applies to  (\ref{eq2}) and the proof for  (\ref{eq3}) applies to  (\ref{eq4}). \eqskip \eqskip \eqskip \hfill $\Box$
\item Consider the case of condition ( B2):  Let $v_{ij}(x)$ be defined as in Theorem \ref{th:theo4.2} -2;

\begin{itemize}[label=$\star$, leftmargin=*]   \setlength\itemsep{1em}
    \item A straightforward calculus shows that by definition of $v_{ij}(x)$ we have:   $r v_{11}(x)-\mathcal{L}_{11} v_{11}(x) - f_{11}(x)= \frac{ K_{21} -K_{11}}{K_{11}} x^{\gamma} - r c_{12}.$
The right hand side of the above equation is positive, hence $v_{ij}(x)$ satisfies (A21). Furthermore, we have: $v_{11}(x)-\big (v_{21}(x) - c_{12} \big )= \hat{V}_{21}(x) - c_{12} -\big (\hat{V}_{21}(x) - c_{12} \big )=0,$ thus, $v_{ij}(x)$ satisfies (A22). With regards to (A23), we have:  $v_{11}(x)- \big (v_{12}(x) +\chi_{12} \big )= \hat{V}_{21}(x) - c_{12} - \big (\hat{V}_{12}(x) - c_{12} +\chi_{12} \big )=-\chi_{12}.$ By Lemma \ref{lem:lemma3}, we deduce that $v_{ij}(x)$ solves (\ref{eq1}).
\item   We have by definition of  $v_{21}(x)$: 
  $r v_{21}(x)-\mathcal{L}_{21} v_{21}(x) - f_{21}(x)=0,$
 hence $v_{ij}(x)$ satisfies (A11). On the other hand, we have:  $v_{21}(x)- \big (v_{11}(x)- c_{12} \big )= c_{12} >0,$ thus, $v_{ij}(x)$ satisfies (A12). With regards to (A13) we have:  $v_{21}(x)- \big (v_{22}(x) + \chi_{21} \big )=- \chi_{21}<0$.  By Lemma \ref{lem:lemma3}, we deduce that $v_{ij}(x)$ solves (\ref{eq3}).
\end{itemize}
We can notice that $v_{11}(x)=v_{12}(x)$ and $v_{21}(x)=v_{22}(x)$, hence the proof for  (\ref{eq1}) applies to  (\ref{eq2}) and the proof for  (\ref{eq3}) applies to  (\ref{eq4}). \hfill $\Box$
\item Consider the case of condition (B3):  Let $v_{ij}(x)$ be defined as in Theorem \ref{th:theo4.2} -3;
\\let us first check solutions for $x <  x^{*} $:
\begin{itemize}[label=$\star$, leftmargin=*]   \setlength\itemsep{1em}
    \item We have by definition  $v_{11}(x)$:   $r v_{11}(x)-\mathcal{L}_{11} v_{11}(x) - f_{11}(x)=r\chi_{12} <0,$ therefore, $v_{ij}(x)$ satisfies (A32). On the other hand, we have: $v_{11}(x)-\big (v_{12}(x)- \chi_{12} \big )=0,$ thus, $v_{ij}(x)$ satisfies (A31). By Lemma \ref{lem:lemma3}, we deduce that $v_{ij}(x)$ solves (\ref{eq1}). 
     \item  We have by definition  $v_{12}(x)$:   $r v_{12}(x)-\mathcal{L}_{12} v_{12}(x) - f_{12}(x)=r\chi_{12} <0,$ hence, $v_{ij}(x)$ satisfies (A12). On the other hand, let $y(x)$ be defined as follows:  $y(x)= v_{12}(x)-\big (v_{22}(x) - c_{12} \big )= ( K_{12}- K_{22})x^{\gamma} + Ax^{m^+_{12}}+ c_{12}.$  A straightforward calculus shows that $y'(x)$ is negative, hence $y(x)$ is decreasing; given that $y(x^*)=0$,  $y$ is positive. Subsequently,  $v_{ij}(x)$ satisfies (A12). With  regards to (A13), we have: $v_{12}(x)- \big (v_{11}(x)+ \chi_{21} \big )=-(\chi_{12} + \chi_{21})<0$.  By Lemma \ref{lem:lemma3}, we deduce that $v_{ij}(x)$ solves (\ref{eq2}).
    \item We have by definition of $v_{ij}(x)$: $r v_{21}(x)-\mathcal{L}_{21} v_{21}(x) - f_{21}(x)=r\chi_{12} < 0,$ therefore, $v_{ij}(x)$ satisfies (A31). On the other hand, we have: $v_{21}(x)-\big (v_{22}(x) +\chi_{12} \big )=0,$ thus, $v_{ij}(x)$ satisfies (A32). By Lemma \ref{lem:lemma3}, we deduce that $v_{ij}(x)$ solves (\ref{eq3}).
    \item We have by definition of $v_{ij}(x)$: $r v_{22}(x)-\mathcal{L}_{22} v_{22}(x) - f_{22}(x)=0,$ therefore, $v_{ij}(x)$ satisfies (A11). Let $s(x)$ be defined as follows: $s(x)=v_{22}(x)-\big (v_{12}(x) -c_{21} \big ) = ( K_{22}- K_{12})x^{\gamma} - Ax^{m^+_{12}}+ c_{21},$  A straightforward calculus shows that $ s'(x)$ is positive and subsequently, $s(x)$ is increasing; given that $s(0)= c_{21}$, $v_{ij}(x)$ satisfies (A12). With regards to (A13), we have: $v_{22}(x)-\big (v_{21}(x) +\chi_{21} \big )=-\chi_{21}<0.$
    By Lemma \ref{lem:lemma3}, we deduce that $v_{ij}(x)$ solves (\ref{eq4}).
\end{itemize}

Let us now check solution for $x \geq  x^{*} $:
\begin{itemize}[label=$\star$, leftmargin=*]   \setlength\itemsep{1em}
    \item  We have by definition of $v_{ij}(x)$,  $ v_{11}(x)-\big (v_{21}(x) - c_{12} \big )=0,$ hence $v_{ij}(x)$ satisfies (A32). On the other hand, $v_{11}(x)-\big (v_{12}(x) + \chi_{12} \big )=\hat{V}_{11}(x) - c_{12} + \chi_{12} - \big ( \hat{V}_{12}(x)  - c_{12} +\chi_{12} \big )=0,$ hence $v_{ij}$ satisfies (A32).  By Lemma \ref{lem:lemma3}, we deduce that $v_{ij}(x)$ solves (\ref{eq1}).
    \item  We have by definition of $v_{ij}(x)$: $r v_{12}(x)-\mathcal{L}_{12} v_{12}(x) - f_{12}(x)=\frac{ K_{22}- K_{12}}{K_{12}} x^{\gamma} -rc_{12},$ A straightforward calculus shows that this quantity is positive, therefore, $v_{ij}(x)$ satisfies (A21). On the other hand, we have: $v_{12}(x)-\big (v_{22}(x) - c_{12} \big )=0,$ thus, $v_{ij}(x)$ satisfies (A22). With regards to (A23), we have: $v_{12}(x)-\big (v_{11}(x) + \chi_{21} \big )=-(\chi_{12} +\chi_{21})<0$. By Lemma \ref{lem:lemma3}, we deduce that $v_{ij}(x)$ solves (\ref{eq2}).
\item We have by definition of $v_{ij}(x)$:  $r v_{21}(x)-\mathcal{L}_{21} v_{21}(x) - f_{21}(x)=r\chi_{12} < 0,$ hence $v_{ij}(x)$ satisfies (A32). On the other hand, we have: $v_{21}(x)-\big (v_{22}(x) +\chi_{12} \big )=0,$ hence  $v_{ij}(x)$ satisfies (A31). By Lemma \ref{lem:lemma3}, we deduce that $v_{ij}(x)$ solves (\ref{eq3}).  
\item  We have by definition of $v_{ij}(x)$: $r v_{22}(x)-\mathcal{L}_{22} v_{22}(x) - f_{22}(x)=0,$ thus, $v_{ij}(x)$ satisfies (A11). On the other hand, we have: $v_{22}(x)-\big (v_{12}(x) - c_{21} \big )= c_{12} + c_{21}>0,$ therefore, $v_{ij}(x)$ satisfies (A12). With regards to (A13), we have: $v_{22}(x)-\big (v_{21}(x) + \chi_{21} \big )= -\big (\chi_{12} + \chi_{21} \big )< 0$.  By Lemma \ref{lem:lemma3}, we deduce that $v_{ij}(x)$ solves (\ref{eq4}).  \hfill $\Box$
\end{itemize}
\item Consider the case of condition (B4):  Let $v_{ij}(x)$ be defined as in Theorem \ref{th:theo4.2} -4: 
\begin{itemize}[label=$\star$, leftmargin=*]   \setlength\itemsep{1em}
    \item We have by definition of $v_{ij}(x)$: 
    $v_{11}(x)-\big (v_{21}(x) - c_{12} \big )= 0,$ thus, $v_{ij}(x)$ satisfies (A32). With regards to (A31), we have: $v_{11}(x)-\big (v_{12}(x) +\chi_{12} \big )=\hat{V}_{11}(x) -c_{12} + \chi_{12}-\big (\hat{V}_{12}(x) - c_{12}+\chi_{12} \big )=0$.  By Lemma \ref{lem:lemma3}, we deduce that $v_{ij}(x)$ solves (\ref{eq1}).
    \item We have by definition of $v_{ij}(x)$: $r v_{12}(x)-\mathcal{L}_{12} v_{12}(x) - f_{12}(x)= \frac{ K_{21} -K_{12}}{K_{12}} x^{\gamma} - r c_{12}.$ The right hand side of the above equation is clearly positive, hence $v_{ij}(x)$ satisfies (A21). Furthermore, we have: $v_{12}(x)-\big (v_{22}(x)-c_{12} \big ) = 0,$ hence (A22) is satisfied. With regards to (A23), we have: $v_{12}(x)-\big (v_{11}(x)+ \chi_{21} \big )= - \chi_{21} -\chi_{12}<0$.  By Lemma \ref{lem:lemma3}, we deduce that $v_{ij}(x)$ solves (\ref{eq2}).
    \item We have by definition of $v_{ij}(x)$: $v_{21}(x)-\big (v_{22}(x)+ \chi_{12} \big )= 0,$ thus $v_{ij}(x)$ satisfies (31). With regards to (A32), we have: $r v_{21}(x)-\mathcal{L}_{21} v_{21}(x) - f_{21}(x)=r\chi_{12} < 0$. By Lemma \ref{lem:lemma3}, we deduce that $v_{ij}(x)$ solves (\ref{eq3}).
    \item We have by definition of $v_{ij}(x)$: $r v_{22}(x)-\mathcal{L}_{22} v_{22}(x) - f_{22}(x)= 0,$ thus,   $v_{ij}(x)$ satisfies (A11). On the other hand, we have: $v_{22}(x)-\big (v_{12}(x) - c_{21} \big )= c_{21} + c_{12} > 0,$ hence $v_{ij}(x)$ satisfies (A12). With regards to (A13), we have : $v_{22}(x)- \big (v_{21}(x) + \chi_{21} \big )= - \chi_{21} -\chi_{12} < 0.$ By Lemma \ref{lem:lemma3}, we deduce that $v_{ij}(x)$ solves (\ref{eq4}).  \hfill $\Box $
\end{itemize}
\end{enumerate}        
\epf
\bethe
\normalfont
\label{th:theo4.3}
 \textit{Case where $K_{11} = K_{12}  >  K_{21} = K_{22}.$ }
\begin{enumerate}[leftmargin=*]
    \item If Condition B1 is verified, then let $v_{ij}(x)$ be defined as follows:
\begin{tcolorbox}[colback=yellow!10!white, colframe=red!50!black]
\begin{minipage}{0.48\textwidth}
\begin{fleqn}[\parindent]
\begin{spreadlines}{1ex}
     \begin{align}
     \; v_{11}(x) & =  \hat{V}_{11}(x), \qquad \qquad \qquad x > 0,\\
    v_{21}(x) &= \begin{cases}
    \hat{V}_{21}(x) +  Ax^{m^+_{21}},& x < x^*\\
    \hat{V}_{11}(x)  -  c_{21}, & x\geq x^* 
\end{cases},
     \end{align}
     \end{spreadlines}
     \end{fleqn}
\end{minipage}
\hfill\vline\hfill
\begin{minipage}{0.5\textwidth}
\begin{fleqn}[\parindent]
\begin{spreadlines}{1ex}
   \begin{align}
   \; v_{12}(x)& =  \hat{V}_{12}(x),\qquad \qquad \qquad x > 0,\\
\; v_{22}(x)& =\begin{cases}
    \hat{V}_{22}(x) +  Ax^{m^+_{22}}, & x < x^*\\
\hat{V}_{11}(x)  -  c_{21}, & x\geq x^*
\end{cases},
     \end{align}
     \end{spreadlines}
     \end{fleqn}
\end{minipage}
\end{tcolorbox}

 where $x^*=\inf S^{\xi}_{21} \in (0, +\infty)$. Using the fact that $v$ is $\mathcal{C}^1$ on $\partial S^{\xi}_{12}$ we get :
$$(K_{11}-K_{21})(x^*)^{\gamma}=\frac{m_{21}^{+} c_{21}}{m_{21}^{+}- \gamma}, \quad
A=(K_{11}-K_{21})\frac{\gamma}{m_{21}^{+}}(x^*)^{\gamma -m_{21}^{+}}.$$
Such $v_{ij}(x)$ are the solutions associated to the system of quasi-variational inequalities with:
    $$S^{\xi}_{11}=S^{\eta}_{11}=S^{\xi}_{12}= S^{\eta}_{12}=S^{\eta}_{21}=S^{\eta}_{22}=\emptyset, \quad S^{\xi}_{21}=S^{\xi}_{22}=[x^*,\infty).$$  
 When the switching costs are positive (condition B1), player I has to stay in regime $2$ for $x$ less than the threshold value $x^*$ and switch beyond that value. Player II on his side has to never switch.
\item If Condition B2 is verified, then let $v_{ij}(x)$ be defined as follows:
\begin{adjustwidth}{-1 cm}{-1 cm}
\begin{tcolorbox}[colback=yellow!10!white, colframe=red!50!black]
\begin{minipage}{0.48\textwidth}
\begin{fleqn}[\parindent]
\begin{spreadlines}{1ex}
     \begin{align}
\;  v_{11}(x)& = \begin{cases}
    \hat{V}_{11}(x) +  Ax^{m^-_{11}}, & x > x^*_A\\
v_{21}(x)  -  c_{12},  & x\leq x^*_A
\end{cases},\\
\; v_{12}(x)& = \begin{cases}
    \hat{V}_{12}(x) +  Ax^{m^-_{12}}, & x > x^*_A\\
v_{22}(x)  -  c_{12},  & x\leq x^*_A
\end{cases},
     \end{align}
     \end{spreadlines}
     \end{fleqn}
\end{minipage}
\hfill\vline\hfill
\begin{minipage}{0.5\textwidth}
\begin{fleqn}[\parindent]
\begin{spreadlines}{1ex}
   \begin{align}
\; v_{21}(x)& = \begin{cases}
    \hat{V}_{21}(x) +  Bx^{m^+_{21}}, & x < x^*_B\\
v_{11}(x)  -  c_{21},  & x \geq x^*_B
\end{cases},\\
\; v_{22}(x)&=\begin{cases}
  \hat{V}_{22}(x) +  Bx^{m^+_{22}},&   x < x^*_B\\
v_{12}(x)  -  c_{21}, &  x \geq x^*_B 
\end{cases},
     \end{align}
     \end{spreadlines}
     \end{fleqn}
\end{minipage}
\end{tcolorbox}
\end{adjustwidth}
 where $x^*_A=\inf S^{\xi}_{12} \in (0, +\infty)$ and $x^*_B=\inf S^{\xi}_{21} \in (0, +\infty)$. Using the fact that $v$ is $\mathcal{C}^1$ on $\partial S^{\xi}_{12}$ we get:
$$\hat{V}_{11}(x^*_A)  + A{x^*_A}^{m^-_{11}}=v_{21}({x^*_A})  - c_{12} = \hat{V}_{21}(x^*_A)  +  B{x^*_A}^{m^+_{21}} - c_{12},$$
$$\hat{V'}_{11}(x^*_A)  + A m^-_{11}{x^*_A}^{m^-_{11}-1}=v'_{21}({x^*_A}) = \hat{V'}_{21}({x^*_A}) + m^+_{21} B{x^*_A}^{m^+_{21}-1},$$
$$\hat{V}_{11}(x^*_B)   + A{x^*_B}^{m^-_{11}}- c_{21}=v_{11}({x^*_B})  - c_{21} = \hat{V}_{21}(x^*_B)   +  B{x^*_B}^{m^+_{21}}, $$
$$\hat{V'}_{11}(x^*_B)   + Am^-_{11}{x^*_B}^{m^-_{11}-1}=v'_{11}({x^*_B})  = \hat{V'}_{21}(x^*_B)   +  m^+_{21}B{x^*_B}^{m^+_{21}-1}, $$
where:
$$ x^*_A= \big[  \frac{-m^-_{11}(c_{21}+c_{12}) \lambda^{m^-_{11}}}{(K_{21} -K_{11})(\gamma - m^+_{21})(1-  \lambda^{m^-_{11}- \gamma} )} \big ] \text{and }  x^*_B=\frac{  x^*_A}{ y},$$
$$A=\frac{(K_{11}-K_{21})(m^+_{11}-\gamma ){x_B^*}^{\gamma-m^-_{21}} + m^+_{11}c_{12}{x_B^*}^{-m^-_{21}} }{m^+_{11}-m^-_{21}},$$
$$ B=A{x_B^*}^{m^-_{21}-m^+_{11}}-(K_{11}-K_{21}){x_B^*}^{\gamma -m^+_{11}}-c_{12}{x_B^*}^{-m^+_{11}}.$$
with $\lambda$ solution in  $\big ( 0,  \big (\frac{-c_{21}}{c_{12}} \big )^{\frac{1}{m^+_{21}}} \big )$ to the equation 
$$m^+_{21}(\gamma -m^-_{11})(1-y^{m^+_{21} - \gamma}) (c_{12}y^{m^-_{11}} + c_{21}) + m^-_{11}(m^+_{21} - \gamma )  \big ( 1-y^{m^-_{11} - \gamma})(c_{12}y^{m^+_{21}} + c_{21})=0.$$

 Such $v_{ij}(x)$ are the solutions associated to the system of quasi-variational inequalities with:
      $$ S^{\xi}_{11}=S^{\xi}_{12}=(0, x^*_A], \quad
 S^{\eta}_{11}= S^{\eta}_{12}=S^{\eta}_{21}= S^{\eta}_{22}=\emptyset, \quad S^{\xi}_{21}=S^{\xi}_{22}=[x^*_B,\infty).$$
When cost of switching from regime $1$ to regime $2$ is negative for player I (condition B2), he has to switch from regime $1$ to regime $2$ when the state variable $x$ is less than the threshold value $x_A^*$, and switch from regime $2$ to regime $1$ if the  state variable $x$ is beyond the threshold value $x_B^*$. Player II on his side has to never switch.
    \item If Condition B3 is verified, then let $v_{ij}(x)$ be defined as follows:

\begin{adjustwidth}{-1 cm}{-1 cm}
\begin{tcolorbox}[colback=yellow!10!white, colframe=red!50!black]
\begin{minipage}{0.5\textwidth}
\begin{fleqn}[\parindent]
\begin{spreadlines}{1ex}
     \begin{align}
    \;  v_{11}(x)& =  \hat{V}_{12}(x) +\chi_{12}, \qquad \qquad \qquad x >0, \\ 
  \; v_{21}(x)& =\begin{cases}
    \hat{V}_{22}(x) +  Ax^{m^+_{22}}+\chi_{12}, & x < x^*\\
\hat{V}_{11}(x)  -  c_{21} + \chi_{12}, & x\geq x^*
\end{cases},
     \end{align}
     \end{spreadlines}
     \end{fleqn}
\end{minipage}
\hfill\vline\hfill
\begin{minipage}{0.48\textwidth}
\begin{fleqn}[\parindent]
\begin{spreadlines}{1ex}
   \begin{align}
   \;  v_{12}(x)&=  \hat{V}_{12}(x), \qquad \qquad \qquad  x >0,\\
 \; v_{22}(x)& = \begin{cases}
    \hat{V}_{22}(x) +  Ax^{m^+_{21}},&  x < x^*\\
\hat{V}_{11}(x)  -  c_{21}, &     x\geq x^*\\
\end{cases},
     \end{align}
     \end{spreadlines}
     \end{fleqn}
\end{minipage}
\end{tcolorbox}
\end{adjustwidth}
where $x^*=\inf S^{\xi}_{21} \in (0, +\infty)$. Using the fact that $v$ is $\mathcal{C}^1$ on $\partial S^{\xi}_{12}$ we get :
  $$(K_{11}-K_{21})(x^*)^{\gamma}=\frac{m_{21}^{+} c_{21}}{m_{21}^{+}- \gamma}, \quad
A=(K_{11}-K_{21})\frac{\gamma}{m_{21}^{+}}(x^*)^{\gamma -m_{21}^{+}}.$$
Such $v_{ij}$ are the solutions associated to the system of quasi-variational inequalities with:
\begin{center}
$S^{\xi}_{11}=S^{\xi}_{12}=S^{\eta}_{12}=S^{\eta}_{22}=\emptyset, \quad
 S^{\eta}_{11}=S^{\eta}_{21}=(0,+\infty), \quad S^{\xi}_{21}=S^{\xi}_{22}=[x^*,+\infty) .$
 \end{center}
When the switching costs are negative from regime $1$ to regime $2$ for player II (condition B3), he has to always switch from regime 1 to regime 2. player I on his side has to switch from regime $2$ to regime $1$ when the state variable $x$ beyond the threshold value $x^*$.
  \item If Condition B4 is verified, then let $v_{ij}(x)$ be defined as follows:

\begin{adjustwidth}{-1 cm}{-1 cm}
\begin{tcolorbox}[colback=yellow!10!white, colframe=red!50!black]
\begin{minipage}{0.53\textwidth}
\begin{fleqn}[\parindent]
\begin{spreadlines}{1ex}
\begin{align}
 v_{11}(x)& = \begin{cases}
    \hat{V}_{11}(x) +  Ax^{m^-_{11}} + \chi_{12},&  x > x^*_A\\
v_{21}(x)  -  c_{12} + \chi_{12}, & x \leq x^*_A
\end{cases},\\
v_{21}(x)& = \begin{cases}
    \hat{V}_{21}(x) +  Bx^{m^+_{21}} + \chi_{12}, & x < x^*_B\\
v_{12}(x)  - c_{21} + \chi_{12},  & x \geq x^*_B
\end{cases},
     \end{align}
     \end{spreadlines}
     \end{fleqn}
\end{minipage}
\hfill\vline\hfill
\begin{minipage}{0.46\textwidth}
\begin{fleqn}[\parindent]
\begin{spreadlines}{1ex}
   \begin{align}
 \;  v_{12}(x)& = \begin{cases}
    \hat{V}_{12}(x) +  Ax^{m^-_{12}},&  x > x^*_A\\
v_{22}(x)  - c_{12},&  x \leq x^*_A
\end{cases},\\
\; v_{22}(x) &= \begin{cases}
    \hat{V}_{22}(x) +  Bx^{m^+_{22}},& x < x^*_A\\
v_{12}(x)  -  c_{21},& x \geq x^*_B
\end{cases},
     \end{align}
     \end{spreadlines}
     \end{fleqn}
\end{minipage}
\end{tcolorbox}
\end{adjustwidth}
 where $x^*_A=\inf S^{\xi}_{12} \in (0, +\infty)$ and $x^*_B=\inf S^{\xi}_{21} \in (0, +\infty)$. Using the fact that $v$ is $\mathcal{C}^1$ on $\partial S^{\xi}_{12}$ we get :
$$\hat{V}_{11}(x^*_A)  + A{x^*_A}^{m^-_{11}}=v_{21}({x^*_A})  - c_{12} = \hat{V}_{21}(x^*_A)  +  B{x^*_A}^{m^+_{21}} - c_{12},$$
$$\hat{V'}_{11}(x^*_A)  + A m^-_{11}{x^*_A}^{m^-_{11}-1}=v'_{21}({x^*_A}) = \hat{V'}_{21}({x^*_A}) + m^+_{21} B{x^*_A}^{m^+_{21}-1},$$
$$\hat{V}_{11}(x^*_B)   + A{x^*_B}^{m^-_{11}}- c_{21}=v_{11}({x^*_B})  - c_{21} = \hat{V}_{21}(x^*_B)   +  B{x^*_B}^{m^+_{21}}, $$
$$\hat{V'}_{11}(x^*_B)   + Am^-_{11}{x^*_B}^{m^-_{11}-1}=v'_{11}({x^*_B})  = \hat{V'}_{21}(x^*_B)   +  m^+_{21}B{x^*_B}^{m^+_{21}-1}, $$
where:
$$ x^*_A= \big[  \frac{-m^-_{11}(c_{21}+c_{12}) \lambda^{m^-_{11}}}{(K_{21} -K_{11})(\gamma - m^+_{21})(1-  \lambda^{m^-_{11}- \gamma} )} \big ] \text{and }  x^*_B=\frac{  x^*_A}{y},$$

with $\lambda$ solution in  $\big ( 0,  \big (\frac{-c_{21}}{c_{12}} \big )^{\frac{1}{m^+_{21}}} \big )$ to the equation 
$$m^+_{21}(\gamma -m^-_{11})(1-y^{m^+_{21} - \gamma}) (c_{12}y^{m^-_{11}} + c_{21}) + m^-_{11}(m^+_{21} - \gamma )  \big ( 1-y^{m^-_{11} - \gamma})(c_{12}y^{m^+_{21}} + c_{21})=0.$$

Such $v_{ij}$ are the solutions associated to the system of quasi-variational inequalities with:
\begin{center}
$S^{\xi}_{11}=S^{\xi}_{12}=(0, x^*_A],\quad
 S^{\eta}_{11}=S^{\eta}_{21}=(0,+\infty), \quad S^{\eta}_{12}=S^{\eta}_{22}=\emptyset, \quad S^{\xi}_{21}= S^{\xi}_{22}=[x^*_B,\infty) .$
\end{center}
When cost of switching from regime $1$ to regime $2$ is negative for both players (condition B4) , player I has to switch from regime $1$ to regime $2$ when the state variable $x$ is less than the threshold value $x_A^*$, and switch from regime $2$ to regime $1$ if the  state variable $x$ is beyond the threshold  $x_B^*$. Player II on his side has to always switch from regime $1$ to regime $2$.
\end{enumerate}
\eethe

\bpf
\normalfont
\begin{enumerate}[leftmargin=*]
    \item  Consider the case of condition (B1): Let $v_{ij}(x)$ be defined as in Theorem \ref{th:theo4.3} -1: 
Let us first check solutions for $x <  x^{*} $:
\begin{itemize}[label=$\star$, leftmargin=*]   \setlength\itemsep{1em}
    \item We have by definition of $v_{ij}(x)$:
    $r v_{11}(x)-\mathcal{L}_{11} v_{11}(x) - f_{11}(x)=0,$
    therefore, $v_{ij}(x)$ satisfies (A11). On the other hand, let $y$ be a function given by $y(x)=  \hat{V}_{11}(x) -\big ( \hat{V}_{21}(x)+Ax^{m_{21}^{+}}-c_{12} \big )= (K_{11}- K_{21})x^{\gamma} - A x^{m_{21}^{+}} + c_{12}.$
A straightforward calculus shows that $y'(x)$ is positive and hence $y$ is increasing. Moreover, 
$y(0)= c_{12}\leq 0$  subsequently $y(x)  > 0$. This implies that $v_{ij}(x)$ satisfies (A12).
With regards to (A13), we have: ${v}_{11}(x)-\big ({v}_{12}(x) + \chi_{12} \big )=  -\chi_{12} < 0.$
By Lemma \ref{lem:lemma3}, we deduce that $v_{ij}(x)$ solves (\ref{eq1}).
\item We have by definition of $v_{ij}(x)$:
    $r v_{21}(x)-\mathcal{L}_{21} v_{21}(x) - f_{21}(x)=0, $
    therefore, $v_{ij}(x)$ satisfies (A11). On the other hand, let $g$ be the function defined as follows: $g(x)= \hat{V}_{21}(x)+Ax^{m_{21}^{+}}-\hat{V}_{11}(x) + c_{21}.$
We can see that $g(x)= -y(x)+c_{21}+c_{12}$,
where $y$ is the function from above, it is clear that $g(x)$ is a decreasing function. Furthermore, $g(x^*) =0$, and therefore, $g$ is a positive function on $[0,x^{*} [$. This implies that $v_{ij}(x)$ satisfies (A12). With regards to (A13), we have: ${v}_{21}(x)-\big ({v}_{22}(x) + \chi_{12} \big )=  -\chi_{12} < 0$. 
By Lemma \ref{lem:lemma3}, we deduce that $v_{ij}(x)$ solves (\ref{eq2}).
\end{itemize}
We can notice that $v_{11}(x)=v_{12}(x)$ and $v_{21}(x)=v_{22}(x)$, so the proof for (\ref{eq1}) applies to (\ref{eq2}) and the proof for  (\ref{eq3}) applies to  (\ref{eq4}).
Let's now check solutions for $x \geq x^{*}$:
\begin{itemize}[label=$\star$, leftmargin=*]   \setlength\itemsep{1em}
    \item We have by definition of $v_{ij}(x)$:
    $r v_{11}(x)-\mathcal{L}_{11} v_{11}(x) - f_{11}(x)=0, $
    thus, $v_{ij}(x)$ satisfies (A11). Furthermore, we have: ${v}_{11}(x)-\big ({v}_{21}(x)-c_{12} \big )= \hat{V}_{11}(x)-\big (\hat{V}_{11}(x)-c_{12} \big )=c_{12}>0,$ hence, $v_{ij}(x)$ satisfies (A12). With regards to (A13), we have: ${v}_{11}(x)-\big ({v}_{12}(x) + \chi_{12} \big )=  -\chi_{12} < 0.$ By Lemma \ref{lem:lemma3}, we deduce that $v_{ij}(x)$ solves (\ref{eq1}).
    \item We have by definition of $v_{ij}(x)$:  $r {v}_{21}(x)-\mathcal{L}_{21} {v}_{21}(x) - f_{21}(x)= \frac{ K_{11}- K_{21}}{K_{21}} x^{\gamma} -rc_{21}.$
Let $h$ be the function defined as follows: $h(x)= \frac{ K_{11}- K_{21}}{K_{21}} x^{\gamma} -rc_{21}.$  $h$ is increasing on $[x^{*}  + \infty)$. 
Furthermore, $h(x^*)= \frac{m_{21}^{+}}{K_{21}(m_{21}^{+}-1)}c_{21}  -rc_{21}.$  A straightforward calculus shows that $\frac{m_{21}^{+}}{K_{21}(m_{21}^{+}-1)}> r$, which implies that $h(x^*)>0$ and therefore $h(x)>0$. This implies that $v_{ij}(x)$ satisfies (A21).  On the other hand, we have:
 ${v}_{21}(x)-\big ({v}_{11}(x)-c_{21} \big )= \hat{V}_{11}(x)- c_{21} -\big (\hat{V}_{11}(x)- c_{21} \big ) =0.$ Therefore, $v_{ij}(x)$ satisfies (A22). With regards to (A23), we have: 
  ${v}_{21}(x)-\big ({v}_{22}(x) + \chi_{12} \big )=  -\chi_{12} < 0.$ By Lemma \ref{lem:lemma3}, we deduce that $v_{ij}(x)$ solves (\ref{eq3}).
\end{itemize}
We can notice that $v_{11}(x)=v_{12}(x)$ and $v_{21}(x)=v_{22}(x)$, so the proof for (\ref{eq1}) applies to (\ref{eq2}) and the proof for (\ref{eq3}) applies to (\ref{eq4}).  \hfill $\Box$
    \item Consider the case of condition (B2):  Let $v_{ij}(x)$ be defined as in Theorem \ref{th:theo4.3}-2:

A straightforward calculus shows that $y<\big (\frac{-c_{21}}{c_{12}} \big )^{\frac{1}{m^+_{11}}} <1$. Hence, $x_A^*<x_B^*$. Now let us show that $v_{ij}(xx)$ defined above solves the system of quasi-variational inequalities:
we shall work respectively on $(0,x_B^*)$,$(x_B^*,+\infty)$. The proof is similar for both ${v}_{11}(x)$ and ${v}_{21}(x)$, so we choose to only prove the result for ${v}_{21}(x)$.

\begin{itemize}[label=$\star$, leftmargin=*]   \setlength\itemsep{1em}
\item By definition, $v_{21}(x)$ satisfies (A11). With regards to (A12),  we have for $x < x_A^*$: $v_{21}(x) -\big (v_{11}(x)-c_{21} \big ) =  v_{21}(x) - v_{21}(x) + c_{21} + c_{12} > 0.$ Now on $(x_A^*,x_B^*)$ we need to prove that: $p(x)=\hat{V}_{21}(x) +  Bx^{m^+_{21}} -  \hat{V}_{11}(x) -  Ax^{m^-_{11}} + c_{21} \geq 0$. We have $p(x_A^*)=c_{21} + c_{12} > 0$ and $p(x_B^*)=0$. Suppose that there exist some $\alpha \in (x_A^*, x_B^*)$  such that $p(\alpha )  =0$; this imply that there exist $\beta \in (\alpha , x_B^*)$ such that $p(\beta )  =0$ As such , the equation $p'(x)=0$ admit at least three solutions in $[x_A^*, x_B^* ]$. However a straightforward study of $p$ shows that $p'$ can take the value zero at most at two points in $(0, +\infty)$. This leads to a contradiction, proving therefore that $p(x) >0$ on $(x_A^*, x_B^*)$; we deduce that $v_{ij}(x)$ satisfies (A12) for $x < x_B^*$. With regards to (A13), we have: $v_{21}(x)-\big (v_{22}(x) + \chi_{12} \big )= - \chi_{12} < 0$. By Lemma \ref{lem:lemma3}, we deduce that $v_{ij}(x)$ solves (\ref{eq3}).
\end{itemize}
Now, for $x \in (x_B^*, +\infty):$

 We will just  show that $v_{ij}(x)$ satisfies (A21) in (\ref{eq3}), since (A22) and (A23) are straightforward.
  \\
  Let $q$ be defined as follows: $q(x)=r {v}_{21}(x)-\mathcal{L}_{21} {v}_{21}(x) - f_{21}(x)=\frac{K_{11}-K_{21}}{K_{21}} x^{\gamma} + m_{11}^- l A x^{m_{11}^-}-rc_{21},$ \\where $l =\frac{1}{2}(\sigma_{11}^2-\sigma_{21}^2)(m_{11}^--1)+b_{11}-b_{21}$.
Let's show that $q(x)>0$ on  $(x_B^*, +\infty)$.
  \begin{itemize}
      \item If $l < 0$, we just have to show that $\frac{K_{11}-K_{21}}{K_{21}} x^{\gamma}-rc_{21}>0$ which is rather straightforward.
      \item If $l\geq 0$, q would be non decreasing on $(0, +\infty)$ with $\lim_{x \rightarrow 0^+}q(x)=-\infty$, and $\lim_{x \rightarrow \infty}q(x)=+\infty$. As such it suffices to show that $q(x_B^*)\geq 0$. 
      $q(x_B^*)=(K_{21}-K_{11}) \big [ \frac{m_{21}^+ - m_{11}^-}{K_{21}}- (m_{21}^+ - \gamma ) m_{11}^-l]-rc_{21}+ m_{21}^+ m_{11}^- c_{21}l$. We can easily check that $q(x_B^*)>\frac{m_{21}^+}{K_{21}(m_{21}^+ -1)} -r>0$, we thus get the desired result.
  \end{itemize}
  
  We can notice that $v_{11}(x)=v_{12}(x)$ and $v_{21}(x)=v_{22}(x)$, so the proof for (\ref{eq1}) applies to (\ref{eq2}) and the proof for (\ref{eq3}) applies to (\ref{eq4}). \eqskip \eqskip \eqskip  \qquad \qquad \hfill $\Box$
    \item Consider the case of condition (B3):  Let $v_{ij}(x)$ be defined as in Theorem \ref{th:theo4.3}-3; let us first check solutions for $x <  x^{*} $:
\begin{itemize}[label=$\star$, leftmargin=*]   \setlength\itemsep{1em}
   \item By definition of $v_{ij}(x)$, we have: $r v_{11}(x)-\mathcal{L}_{11} v_{11}(x) - f_{11}(x)=r\chi_{12} < 0.$ Therefore,  $v_{ij}(x)$ satisfies (A32). On the other hand:
   ${v}_{11}(x)-\big ({v}_{12}(x)-\chi_{12} \big )=0$, hence $v_{ij}(x)$ satisfies (A31). By Lemma \ref{lem:lemma3}, we deduce that $v_{ij}(x)$ solves (\ref{eq1}).
    \item By definition of $v_{ij}(x)$, we have: 
    $r v_{12}(x)-\mathcal{L}_{12} v_{12}(x) - f_{12}(x)=0.$
    Therefore, $v_{ij}(x)$ satisfies (A11). On the other hand,  $v_{12}(x)- \big (v_{22}(x) - c_{12} \big )=\hat{V}_{12}(x) - \big ( \hat{V}_{22}(x)+ Ax^{m^+_{22}}  - c_{12} \big ).$ Let $p$ be as follows: $p(x)=\hat{V}_{12}(x) - \hat{V}_{22}(x) -Ax^{m^+_{22}} + c_{12}$. We can check that we also have 
$p(x)= (K_{12}- K_{22}) x^{\gamma} -Ax^{m^+_{22}} + c_{12}$.
The second derivative of $p$ is negative and therefore $p'(x)$ is decreasing. On the other hand, $p'(x^*)= 0$, therefore $p'(x)$ is a positive function, which implies that $p$ is increasing. Given that $p(0)= c_{12}$,  we deduce that $v_{ij}(x)$ satisfies (A12). With regards to (A13):
$v_{12}(x)- \big (v_{11}(x) + \chi_{21} \big )=\hat{V}_{12}(x)-\big (\hat{V}_{12}(x)+\chi_{12}+\chi_{21} \big )<0.$
By Lemma \ref{lem:lemma3}, we deduce that $v_{ij}(x)$ solves (\ref{eq2}).
    \item By definition of $v_{ij}(x)$, we have:  $v_{21}(x)- \big (v_{22}(x) +\chi_{12} \big )= \hat{V}_{22}(x) + Ax^{m^+_{22}}+\chi_{12}-\big (\hat{V}_{22}(x) + Ax^{m^+_{22}} +\chi_{12} \big )=0.$
 Therefore, $v_{ij}(x)$ satisfies (A31). On the other hand, let $q$ be defined as follows: $q(x)=v_{21}(x)-\big (v_{11}(x)-c_{21} \big ).$
 We can see that $q= -p$, where $p$ is the function from above; Given that $p$ is positive on $(0, x^*)$, $-p$ is negative, therefore $v_{ij}(x)$ satisfies (A32). By Lemma \ref{lem:lemma3}, we deduce that $v_{ij}(x)$ solves (\ref{eq3}).
\item By definition of $v_{ij}(x)$, we have: $r v_{22}(x)-\mathcal{L}_{22} v_{22}(x) - f_{22}(x)=0. $ Therefore, $v_{ij}(x)$ satisfies (A11). On the other hand, we have: $v_{22}(x)-\big (v_{12}(x)-c_{21} \big )=q(x)$ where  $q(x)$ is the function from above; we deduce that $v_{ij}(x)$ satisfies (A12). With regards to (A13), wwe have:  $v_{22}(x)-\big (v_{21}(x)+\chi_{21} \big )=- \big (\chi_{21} +\chi_{12} \big )<0$. By Lemma \ref{lem:lemma3}, we deduce that $v_{ij}(x)$ solves (\ref{eq4}).
\end{itemize}

Let's check solutions for $x \geq  x^{*} $:

\begin{itemize}[label=$\star$, leftmargin=*]   \setlength\itemsep{1em}
    \item Using same steps as for $x < x^*$, we prove that $v_{ij}(x)$ solves (\ref{eq1}).
    \item By definition of $v_{ij}(x)$, we have:  $r v_{12}(x)-\mathcal{L}_{12} v_{12}(x) - f_{12}(x)=0.$ Therefore,  $v_{ij}(x)$ satisfies (A11). On the other hand, $v_{12}(x)- \big (v_{22}(x) - c_{12} \big )=\hat{V}_{12}(x) - \big ( \hat{V}_{22}(x)  - c_{21} -c_{12}\big )=c_{12}+c_{21} > 0, $ therefore $v_{ij}(x)$ satisfies (A12). With regards to (A13): $v_{12}(x)- \big (v_{11}(x) + \chi_{21} \big )= - \big (\chi_{21} + \chi_{12} \big ) <0.$ By Lemma \ref{lem:lemma3}, we deduce that $v_{ij}(x)$ solves (\ref{eq2}).
    \item By definition of $v_{ij}(x)$, we have:  $v_{21}(x)-\big (v_{11}(x)-c_{21} \big )=0,$ therefore $v_{ij}(x)$ satisfies (A31). On the other hand, we have: $v_{21}(x)- \big (v_{22}(x) +\chi_{12} \big )= \hat{V}_{11}(x)  - c_{21} + \chi_{12}-\big (\hat{V}_{11}(x)  -  c_{21} + \chi_{12} \big )=0,$ hence $v_{ij}(x)$ satisfies (A32). By Lemma \ref{lem:lemma3}, we deduce that $v_{ij}(x)$ solves (\ref{eq3}).
    \item By definition of $v_{ij}(x)$, we have: $v_{22}(x)-\big (v_{12}(x)-c_{21} \big )=0,$ thus $v_{ij}(x)$ satisfies (A22). On the other hand, we have:  $rv_{22}(x)-\mathcal{L}_{22} v_{22}(x) - f_{22}(x)= \frac{K_{11} - K_{22}}{K_{22}} x^{\gamma} - rc_{21}.$
Let $u$ be defined as follows: $u(x)=\frac{K_{11} - K_{22}}{K_{22}} x^{\gamma} - rc_{21}.$ A straightforward calculus shows that  $u$ is an increasing function. Moreover,  we have:
$u(x^*)= \big (\frac{m^+_{22}}{K_{22}(m^+_{22}-\gamma ) } -r \big ) c_{21}$.
We can easily see that $\big (\frac{m^+_{22}}{K_{22}(m^+_{22}-\gamma ) } -r \big ) > 0$, which implies that $v_{ij}(x)$ satisfies (A21). With regards to (A23), we have $v_{22}(x)-\big (v_{21}(x)+\chi_{21} \big )=-(\chi_{12}+\chi_{21})<0.$ By Lemma \ref{lem:lemma3}, we deduce that $v_{ij}(x)$ solves (\ref{eq4}). \hfill $\Box$
\end{itemize}
    \item Consider the case of condition (B4):  Let $v_{ij}(x)$ be defined as in Theorem \ref{th:theo4.3}-4\\
The proof is similar  for $ v_{21}(x) $ and $v_{11}(x)$ on the one hand, and for $ v_{12}(x)$ and $v_{22}(x)$ on the other hand. Hence we choose to prove results for $v_{21}(x)$ and $v_{22}(x)$. we shall work respectively on $(0,x_B^*)$,$(x_B^*,+\infty)$.\\
Let's check solution for $x< x_B^*$:
\begin{itemize}[label=$\star$, leftmargin=*]   \setlength\itemsep{1em}
    \item By definition of $v_{ij}(x)$, we have:   $v_{21}(x)-\big (v_{22}(x)+\chi_{12} \big )=0,$ therefore $v_{ij}(x)$ satisfies (A31). On the other hand, we have:
    $rv_{21}(x)-\mathcal{L}_{21} v_{21}(x) - f_{21}(x)=r \chi_{12}< 0,$ hence  $v_{ij}(x)$ satisfies (A32). By Lemma \ref{lem:lemma3}, we deduce that $v_{ij}(x)$ solves (\ref{eq3}).
    \item By definition of $v_{ij}(x)$, we have:  $rv_{22}(x)-\mathcal{L}_{22} v_{22}(x) - f_{22}(x)=0, $ therefore, $v_{ij}(x)$ satisfies (A11). on the other hand, we have:
    $v_{22}(x)-\big (v_{12}(x)-c_{21} \big )= \hat{V}_{21}(x) +  Bx^{m^+_{21}} -  \hat{V}_{11}(x) -  Ax^{m^-_{11}} + c_{21}.$ By using the proof in 4.3-2, we know that this quantity is positive, hence $v_{ij}(x)$ satisfies (A12). With regards to (A13), we have: $v_{22}(x)-\big (v_{21}(x)- \chi_{21} \big )= -\chi_{21} - \chi_{12} < 0.$ By Lemma \ref{lem:lemma3}, we deduce that $v_{ij}(x)$ solves (\ref{eq4}).
\end{itemize}
Let's check solution for $x \geq x_B^*$:
\begin{itemize}[label=$\star$, leftmargin=*]   \setlength\itemsep{1em}
    \item By definition of $v_{ij}(x)$ we have: $v_{21}(x)-\big (v_{22}(x)+ \chi_{12} \big )=0,$ therefore $v_{ij}(x)$ satisfies (A31).  On the other hand,  $v_{21}(x)-\big (v_{11}(x)-c_{21} \big )=0,$ hence $v_{ij}(x)$ satisfies (A32). By Lemma \ref{lem:lemma3}, we deduce that $v_{ij}(x)$ solves (\ref{eq3}).
    \item Using Proof 4.3-2, we have $rv_{22}(x)-\mathcal{L}_{22} v_{22}(x) - f_{22}(x)\geq 0,$ therefore $v_{ij}(x)$ satisfies (A21). On the other hand, we have:  $v_{22}(x)-\big (v_{12}(x)-c_{21} \big )=0,$ therefore, $v_{ij}(x)$ satisfies (A22). With respect to (A23), we have: $v_{22}(x)-\big (v_{21}(x)+ \chi_{21} \big )= - \chi_{21} - \chi_{12} < 0$. By Lemma \ref{lem:lemma3}, we deduce that $v_{ij}(x)$ solves (\ref{eq4}). \hfill $\Box$
\end{itemize}
\end{enumerate}
\epf
\bethe
\normalfont Case where $K_{11} =   K_{21} <  K_{12} = K_{22}.$
\label{th:theo4.4}
\begin{enumerate}[leftmargin=*]
    \item If Condition B1 is verified, then let $v_{ij}(x)$ be defined as follows:

\begin{tcolorbox}[colback=yellow!10!white, colframe=red!50!black]
\begin{minipage}{0.48\textwidth}
\begin{fleqn}[\parindent]
\begin{spreadlines}{1ex}
     \begin{align}
     v_{11}(x)& =\hat{V}_{11}(x), \qquad  \qquad \qquad x >0,\\
     v_{12}(x)& = \begin{cases}
    \hat{V}_{12}(x) +  Ax^{m^+_{12}}, & x <x^*\\
\hat{V}_{11}(x)  +  \chi_{21},  & x \geq x^* 
\end{cases}, 
     \end{align}
     \end{spreadlines}
     \end{fleqn}
\end{minipage}
\hfill\vline\hfill
\begin{minipage}{0.5\textwidth}
\begin{fleqn}[\parindent]
\begin{spreadlines}{1ex}
   \begin{align}
   \; v_{21}(x)& = \hat{V}_{21}(x), \qquad  \qquad \qquad x >0,\\
\; v_{22}(x)& = \begin{cases}
\hat{V}_{22}(x) +  Ax^{m^+_{22}}, & x <x^*\\
\hat{V}_{21}(x)  +  \chi_{21},  & x \geq x^* 
 \end{cases},
     \end{align}
     \end{spreadlines}
     \end{fleqn}
\end{minipage}
\end{tcolorbox}
where $x^*=\inf S^{\eta}_{21} \in (0, +\infty)$. Using the fact that $v$ is $\mathcal{C}^1$ on $\partial S^{\xi}_{12}$ we get:
$$(K_{12}-K_{11})(x^*)^{\gamma}=\frac{m_{12}^{+} \chi_{12}}{m_{12}^{+}- \gamma}, \quad
A=(K_{11}-K_{12})\frac{\gamma}{m_{12}^{+}}(x^*)^{\gamma -m_{12}^{+}}.$$
Such $v_{ij}(x)$ are the solutions associated to the system of quasi-variational inequalities with
$$S^{\xi}_{11}=S^{\eta}_{11}= S^{\xi}_{12}=S^{\eta}_{21}=S^{\xi}_{22}=S^{\xi}_{21}=\emptyset, \quad S^{\eta}_{12}=S^{\eta}_{22}=[x^*,\infty).$$
When the switching costs are positive (condition B1), player II has to stay in regime 2 for x less than the threshold value $x^*$ and switch beyond that value.
\item If Condition B2 is verified then  let $v_{ij}(x)$ be defined as follows:
\begin{adjustwidth}{-0.75 cm}{-0.75 cm}
\begin{tcolorbox}[colback=yellow!10!white, colframe=red!50!black]
\begin{minipage}{0.48\textwidth}
\begin{fleqn}[\parindent]
\begin{spreadlines}{1ex}
     \begin{align}
 v_{11}(x) &=\hat{V}_{21}(x)- c_{12},\qquad \qquad \qquad x > 0,\\
 v_{12}(x)& =\begin{cases}
\hat{V}_{12}(x) +  Ax^{m^+_{12}} - c_{12}, & x <x^*\\
\hat{V}_{11}(x)  -   c_{12}+ \chi_{21},  & x \geq x^* 
\end{cases}, 
     \end{align}
     \end{spreadlines}
     \end{fleqn}
\end{minipage}
\hfill\vline\hfill
\begin{minipage}{0.5\textwidth}
\begin{fleqn}[\parindent]
\begin{spreadlines}{1ex}
   \begin{align}
   \; v_{21}(x)&= \hat{V}_{21}(x), \qquad \qquad \qquad \quad \hspace{0.3 cm} x > 0 ,\\
\; v_{22}(x)&= \begin{cases}
\hat{V}_{22}(x) +  Ax^{m^+_{22}}, & \qquad x <x^*\\
\hat{V}_{21}(x)  +  \chi_{21}, &  \qquad x \geq x^*
 \end{cases},
     \end{align}
     \end{spreadlines}
     \end{fleqn}
\end{minipage}
\end{tcolorbox}
\end{adjustwidth}
where $x^*=\inf S^{\eta}_{21} \in (0, +\infty)$. Using the fact that $v$ is $\mathcal{C}^1$ on $\partial S^{\xi}_{12}$ we get:
$$(K_{12}-K_{11})(x^*)^{\gamma}=\frac{m_{12}^{+} \chi_{12}}{m_{12}^{+}- \gamma}, \quad
A=(K_{11}-K_{12})\frac{\gamma}{m_{12}^{+}}(x^*)^{\gamma -m_{12}^{+}}.$$
Such $v_{ij}(x)$ are the solutions associated to the system of quasi-variational inequalities with: 
$$S^{\xi}_{11}=S^{\xi}_{12}=(0,+\infty), \quad S^{\eta}_{12}=S^{\eta}_{22}=[x^*,+\infty), \quad S^{\xi}_{21}=S^{\eta}_{21}= S^{\xi}_{22}=S^{\eta}_{11}=\emptyset.$$
When the switching cost from regime $1$ to regime $2$ is negative for player I (condition B2), he has to always switch from regime $1$ to regime $2$. Player II on his side has to stay in the continuation for $x$ less than the threshold value $x^*$ and switch beyond that value.
\item If Condition B3 is verified then  let $v_{ij}(x)$ be defined as follows:

\begin{tcolorbox}[colback=yellow!10!white, colframe=red!50!black]
\begin{minipage}{0.48\textwidth}
\begin{fleqn}[\parindent]
\begin{spreadlines}{1ex}
     \begin{align}
    v_{11}(x)& =\begin{cases}
    \hat{V}_{11}(x) +  Ax^{m^-_{11}}, & x >x_A^*\\
v_{12}(x)  +  \chi_{12},  & x \leq x_A^* 
\end{cases},\\
v_{12}(x)&= \begin{cases}
    \hat{V}_{12}(x) +  Bx^{m^+_{12}}, & x <x_B^*\\
v_{11}(x)  +  \chi_{21},  & x \geq x_B^*
\end{cases},
     \end{align}
     \end{spreadlines}
     \end{fleqn}
\end{minipage}
\hfill\vline\hfill
\begin{minipage}{0.5\textwidth}
\begin{fleqn}[\parindent]
\begin{spreadlines}{1ex}
   \begin{align}
   v_{21}(x)& = \begin{cases}
    \hat{V}_{21}(x) +  Ax^{m^-_{21}}, & x >x_A^*\\
v_{22}(x)  +  \chi_{12},  & x \leq x_A^* 
\end{cases},\\
v_{22}(x)&= \begin{cases}
  \hat{V}_{22}(x) +  Bx^{m^+_{22}}, & x <x_B^*\\
v_{21}(x)  +  \chi_{21},  & x \geq x_B^*  
 \end{cases},
     \end{align}
     \end{spreadlines}
     \end{fleqn}
\end{minipage}
\end{tcolorbox}

 where $x^*_A=\inf S^{\eta}_{12} \in (0, +\infty)$ and $x^*_B=\inf S^{\eta}_{21} \in (0, +\infty)$. Using the fact that $v$ is $\mathcal{C}^1$ on $\partial S^{\xi}_{12}$ we get:
$$\hat{V}_{12}(x^*_B)  + B{x^*_B}^{m^+_{12}}=v_{11}({x^*_B})  + \chi_{21} = \hat{V}_{11}(x^*_B)  +  A{x^*_B}^{m^-_{11}} + \chi_{21},$$
$$\hat{V'}_{12}(x^*_B)  + B m^+_{12}{x^*_B}^{m^+_{12}-1}=v'_{11}({x^*_B}) = \hat{V'}_{11}({x^*_B}) + m^-_{11} A{x^*_B}^{m^-_{11}-1},$$
$$\hat{V}_{11}(x^*_A)  + A{x^*_A}^{m^-_{11}}=v_{12}({x^*_A})  + \chi_{12} = \hat{V}_{12}(x^*_A)  +  B{x^*_A}^{m^+_{12}} + \chi_{12},$$
$$\hat{V'}_{11}(x^*_A)  + A m^-_{11}{x^*_A}^{m^-_{11}-1}=v'_{12}({x^*_A}) = \hat{V'}_{12}({x^*_A}) + m^+_{12} B{x^*_A}^{m^-_{12}-1},$$

where:
$$ x^*_A= \big[  \frac{m^-_{11}(\chi_{21}+ \chi_{12}) \lambda^{m^+_{12}}}{(K_{12} -K_{11})(\gamma - m^-_{11})(1-  \lambda^{m^+_{12}- \gamma} )} \big ]^{\frac{1}{\gamma}} \text{and }  x^*_B=\frac{  x^*_A}{ \lambda},$$

with $\lambda$ solution in  $\big ( 0,  \big (\frac{-\chi_{21}}{\chi_{12}} \big )^{\frac{1}{m^+_{12}}} \big )$ to the equation 
$$m^+_{11}(\gamma -m^-_{21})(1-y^{m^+_{11} - \gamma}) (\chi_{12} y^{m^-_{21}} + \chi_{21}) + m^-_{21}(m^+_{11} - \gamma )  \big ( 1-y^{m^-_{21} - \gamma})(\chi_{12} y^{m^+_{11}} + \chi_{21})=0.$$

Such $v_{ij}(x)$ are the solutions associated to the system of quasi- variational inequalities with
\begin{center}
$S^{\xi}_{11}= S^{\xi}_{21}= S^{\xi}_{12}=  S^{\xi}_{22}=\emptyset, \quad
 S^{\eta}_{11}=S^{\eta}_{21}=[0, x^*_A], \quad  S^{\eta}_{12}=S^{\eta}_{22}=[x^*_B,\infty) .$
\end{center}
When cost of switching from regime $1$ to regime $2$ is negative for player II (condition B3), he has to switch from regime $1$ to regime $2$ when the state variable $x$ is less than the threshold value $x_A^*$, and switch from regime $2$ to regime $1$ if the  state variable $x$ is beyond the threshold value $x_B^*$. Player I on his side has to never switch.

\item If Condition B4 is verified then  let $v_{ij}(x)$ be defined as follows:

\begin{adjustwidth}{-0.75 cm}{-0.75 cm}
\begin{tcolorbox}[colback=yellow!10!white, colframe=red!50!black]
\begin{minipage}{0.5\textwidth}
\begin{fleqn}[\parindent]
\begin{spreadlines}{1ex}
     \begin{align}
    v_{11}(x)& =\begin{cases}
    \hat{V}_{11}(x) +  Ax^{m^-_{11}}  -  c_{12}, & x >x_A^*\\
v_{12}(x)  - c_{12}  +  \chi_{12},  & x \leq x_A^* 
\end{cases},\\
 v_{12}(x)& = \begin{cases}
    \hat{V}_{12}(x) +  Bx^{m^+_{12}} - c_{12}, & x <x_B^*\\
v_{11}(x)  + \chi_{21} - c_{12},  & x \geq x_B^* 
\end{cases},
     \end{align}
     \end{spreadlines}
     \end{fleqn}
\end{minipage}
\hfill\vline\hfill
\begin{minipage}{0.48\textwidth}
\begin{fleqn}[\parindent]
\begin{spreadlines}{1ex}
   \begin{align}
\; v_{21}(x)& =\begin{cases}
    \hat{V}_{21}(x) +  Ax^{m^-_{21}}, & x >x_A^*\\
v_{22}(x)  +  \chi_{12},  & x \leq x_A^* 
\end{cases},\\
\; v_{22}(x)& =\begin{cases}
     \hat{V}_{22}(x) +  Bx^{m^+_{22}}, & x <x_B^*\\
v_{21}(x)  +  \chi_{21},  & x \geq x_B^* 
 \end{cases},
     \end{align}
     \end{spreadlines}
     \end{fleqn}
\end{minipage}
\end{tcolorbox}
\end{adjustwidth}

 where $x^*_A=\inf S^{\eta}_{12} \in (0, +\infty)$ and $x^*_B=\inf S^{\eta}_{21} \in (0, +\infty)$. Using the fact that $v$ is $\mathcal{C}^1$ on $\partial S^{\xi}_{12}$ we get :
$$\hat{V}_{12}(x^*_B)  + B{x^*_B}^{m^+_{12}}=v_{11}({x^*_B})  + \chi_{21} = \hat{V}_{11}(x^*_B)  +  A{x^*_B}^{m^-_{11}} + \chi_{21},$$
$$\hat{V'}_{12}(x^*_B)  + B m^+_{12}{x^*_B}^{m^+_{12}-1}=v'_{11}({x^*_B}) = \hat{V'}_{11}({x^*_B}) + m^-_{11} A{x^*_B}^{m^-_{11}-1},$$
$$\hat{V}_{11}(x^*_A)  + A{x^*_A}^{m^-_{11}}=v_{12}({x^*_A})  + \chi_{12} = \hat{V}_{12}(x^*_A)  +  B{x^*_A}^{m^+_{12}} + \chi_{12},$$
$$\hat{V'}_{11}(x^*_A)  + A m^-_{11}{x^*_A}^{m^-_{11}-1}=v'_{12}({x^*_A}) = \hat{V'}_{12}({x^*_A}) + m^+_{12} B{x^*_A}^{m^-_{12}-1},$$

where:
$$ x^*_A= \big[  \frac{m^-_{11}(\chi_{21}+ \chi_{12}) \lambda^{m^+_{12}}}{(K_{12} -K_{11})(\gamma - m^-_{11})(1-  \lambda^{m^+_{12}- \gamma} )} \big ]^{\frac{1}{\gamma}} \text{and }  x^*_B=\frac{ x^*_A}{ \lambda},$$

with $\lambda$ solution in  $\big ( 0,  \big (\frac{-\chi_{21}}{\chi_{12}} \big )^{\frac{1}{m^+_{12}}} \big )$ to the equation
$$m^+_{11}(\gamma -m^-_{21})(1-y^{m^+_{11} - \gamma}) (\chi_{12} y^{m^-_{21}} + \chi_{21}) + m^-_{21}(m^+_{11} - \gamma )  \big ( 1-y^{m^-_{21} - \gamma})(\chi_{12} y^{m^+_{11}} + \chi_{21})=0.$$
Such $v_{ij}(x)$ are the solutions associated to the system of quasi- variational inequalities with
\begin{center}
    \qquad $S^{\xi}_{11}= S^{\xi}_{12}=S^{\eta}_{21}=(0,+\infty),\quad 
 S^{\eta}_{11}= S^{\xi}_{21}=(0, x^*_A], \quad  S^{\eta}_{12}=S^{\eta}_{22}=[x^*_B,\infty) \quad  S^{\xi}_{22}= \emptyset .$  
\end{center}
When cost of switching from regime $1$ to regime $2$ is negative for both players(condition B4) , player II has to switch from regime $1$ to regime $2$ when the state variable $x$ is less than the threshold value $x_A^*$, and switch from regime $2$ to regime $1$ if the  state variable $x$ is beyond the threshold value $x_B^*$. Player I on his side has to always switch from regime $1$ to regime $2$.
\end{enumerate} 
\eethe
 \bpf
 \normalfont
 \begin{enumerate}[leftmargin=*]
     \item Consider the case of Condition (B1):  Let $v_{ij}(x)$ be defined as in Theorem \ref{th:theo4.4}-1 Let us check solutions for $x <  x^{*} $:
\begin{itemize}[label=$\star$, leftmargin=*]   \setlength\itemsep{1em}
    \item By definition of $v_{ij}(x)$, we have:
    $r v_{11}(x)-\mathcal{L}_{11} v_{11}(x) - f_{11}(x)=0,$
    therefore, $v_{ij}(x)$ satisfies (A11). On the other hand,  we have: $v_{11}(x)-\big (v_{21}(x)- c_{12} \big )= c_{12} >0,$
    hence, $v_{ij}(x)$ satisfies (A12). With regards to (A13),  let $z$ be deined as follows: $z(x)=\hat{V}_{11}(x)-\hat{V}_{12}(x)- Ax^{m^+_{12}}-\chi_{12}= (K_{11} -K_{12}) x ^{\gamma}- Ax^{m^+_{12}}-\chi_{12}. $
    $z(x)$ is clearly negative, which implies that $v_{ij}(x)$ satisfies (A13).
    By Lemma \ref{lem:lemma3}, we deduce that $v_{ij}(x)$ solves (\ref{eq1}).
    \item By definition of $v_{ij}(x)$, we have:
    $r v_{12}(x)-\mathcal{L}_{12} v_{12}(x) - f_{12}(x)=0,$  therefore, $v_{ij}(x)$ satisfies (A11). On the other hand,  we have: 
     $v_{12}(x)-\big (v_{22}(x)- c_{12} \big )= c_{12} >0,$ hence $v_{ij}(x)$ satisfies (A12). With regards to (A13), we have: $v_{12}(x)- \big (v_{11}(x) + \chi_{21} \big )=(K_{12} -K_{11}) x ^{\gamma} + Ax^{m^+_{12}} - \chi_{21}.$
Let $y$ be defined as follows: $y(x)=(K_{12} -K_{11}) x ^{\gamma} + Ax^{m^+_{12}} - \chi_{21},$  $y'$ is clearly positive, with $y(x^*)=0$ and hence $v_{ij}$ satisfies (A13). By Lemma \ref{lem:lemma3}, we deduce that $v_{ij}(x)$ solves (\ref{eq2}).\\
 We can notice that $v_{11}(x)=v_{21}(x)$ and $v_{12}(x)=v_{22}(x)$, so the proof for (\ref{eq1}) applies to (\ref{eq3}) and the proof for (\ref{eq2}) applies to (\ref{eq4}).
\end{itemize}
Let's now check solution for $x \geq  x^{*} $:
\begin{itemize}[label=$\star$, leftmargin=*]   \setlength\itemsep{1em}
\item  By definition of $v_{ij}(x)$, we have:
    $r v_{11}(x)-\mathcal{L}_{11} v_{11}(x) - f_{11}(x)=0,$
    therefore, $v_{ij}(x)$ satisfies (A11). On the other hand,  we have: $v_{11}(x)-\big (v_{21}(x)- c_{12} \big )= c_{12} >0,$
    hence, $v_{ij}(x)$ satisfies (A12). With regards to (A13), we have: $v_{11}(x)-\big (v_{12}(x)+ \chi_{12} \big )=-(\chi_{12} +\chi_{21})<0.$ By Lemma \ref{lem:lemma3}, we deduce that $v_{ij}(x)$ solves (\ref{eq1}).
\item By definition of $v_{ij}(x)$, we have: $v_{12}(x)-\big (v_{11}(x)+ \chi_{21} \big )= 0 .$ With regards to (A32), we have: $r v_{12}(x)-\mathcal{L}_{12} v_{12}(x) - f_{12}(x)=\frac{ K_{11} -K_{12}}{K_{12}} x^{\gamma} + r \chi_{21},$ A straightforward calculus shows that this quantity is negative on $[x^{*}, +\infty [$ and we are done. By Lemma \ref{lem:lemma3}, we deduce that $v_{ij}(x)$ solves (\ref{eq2}).
\end{itemize}
Again, given that $v_{11}(x)=v_{21}(x)$ and $v_{12}(x)=v_{22}(x)$, the proof for (\ref{eq1}) applies to (\ref{eq3}) and the proof for (\ref{eq2}) applies to (\ref{eq4}). \hfill $\Box$
\item Consider the case of condition (B2):  Let $v_{ij}(x)$ be defined as in Theorem \ref{th:theo4.4}-2
Let us check solutions for $x <  x^{*} $:

\begin{itemize}[label=$\star$, leftmargin=*]   \setlength\itemsep{1em}
    \item By definition of $v_{ij}(x)$, we have:  $r v_{11}(x)-\mathcal{L}_{11} v_{11}(x) - f_{11}(x)= -rc_{12} >0,$ hence $v_{ij}(x)$ satisfies (A21). On the other hand,  $v_{11}(x)-\big (v_{21}(x) - c_{12} \big )= \hat{V}_{21}(x) - c_{12} -\big (\hat{V}_{21}(x)-c_{12} \big ) =0,$ hence $v_{ij}(x)$ satisfies (A22). With regards to (A23), we have: $v_{11}(x)-\big (v_{12}(x) + \chi_{12} \big )= (K_{11} -K_{12}) x ^{\gamma} - Ax^{m^+_{12}} - \chi_{21};$ 
    let $g$ be defined as follows: $g(x)=(K_{11} -K_{12}) x ^{\gamma} - Ax^{m^+_{12}} - \chi_{12},$  $g(x)$  is clearly negative. By Lemma \ref{lem:lemma3}, we deduce that $v_{ij}(x)$ solves (\ref{eq1}).
    \item By definition of $v_{ij}(x)$, we have: $r v_{12}(x)-\mathcal{L}_{12} v_{12}(x) - f_{12}(x)= -rc_{12} >0,$ therefore $v_{ij}(x)$ satisfies (A21). On the other hand, we have: $v_{12}(x)-\big(v_{22}(x)- c_{12} \big)= 0.$ therefore $v_{ij}(x)$ satisfies (A22).  With regards to (A23), we have:  $v_{12}(x)- (v_{11}(x) + \chi_{21})=(K_{12} -K_{11}) x ^{\gamma} + Ax^{m^+_{12}} - \chi_{21}$.  
    Let $z$ be defined as follows:  $z(x)=(K_{12} -K_{11}) x ^{\gamma} + Ax^{m^+_{12}} - \chi_{21}$. A straightforward calculus shows that $z'(x)>0$, hence $z(x)$ is increasing. Moreover $z(x^*)=0$, hence $z(x)$ is negative and thus $v_{ij}(x)$ satisfies (A23). By Lemma \ref{lem:lemma3}, we deduce that $v_{ij}(x)$ solves (\ref{eq1}).
    \item By definition of $v_{ij}(x)$, we have: $r v_{21}(x)-\mathcal{L}_{21} v_{21}(x) - f_{21}(x)= 0,$ therefore $v_{ij}(x)$ satisfies (A11). On the other hand, we have: $v_{21}(x)- \big (v_{11}(x) -c_{21} \big )=c_{21}+c_{12} >0,$ therefore $v_{ij}(x)$ satisfies (A12). With regards to (A13), we have:  $v_{21}(x) - (v_{22}(x) +\chi_{12})=(K_{21} -K_{22}) x ^{\gamma} - Ax^{m^+_{22}} - \chi_{12},$ The above quantity is clearly negative and hence $v_{ij}(x)$ satisfies (A13). By Lemma \ref{lem:lemma3}, we deduce that $v_{ij}(x)$ solves (\ref{eq3}).
    \item By definition of $v_{ij}(x)$, we have: 
    $r v_{11}(x)-\mathcal{L}_{11} v_{11}(x) - f_{11}(x)=0,$ therefore $v_{ij}(x)$ satisfies (A11). On the other hand, $v_{22}(x)- \big (v_{12}(x) -c_{21} \big )=c_{21}+c_{12} >0,$  hence $v_{ij}(x)$ satisfies (A12). With regars to (A13), we have: $v_{22}(x)- \big(v_{21}(x) + \chi_{21} \big )=(K_{22} -K_{21}) x ^{\gamma} + Ax^{m^+_{22}} - \chi_{21}.$
    Let $p$ be defined as follows: $p(x)=(K_{22} -K_{21}) x ^{\gamma} + Ax^{m^+_{22}} - \chi_{21},$ A straightforward calculus shows that $p'(x)$ is positive, hence $p(x)$ is increasing. Moreover $p(x^*)=0$, hence $p(x)$ is negative $v_{ij}(x)$ satisfies (A13).  By Lemma \ref{lem:lemma3}, we deduce that $v_{ij}(x)$ solves (\ref{eq4}).
\end{itemize}
Let's check solution for $x \geq  x^{*} $:
\begin{itemize}[label=$\star$, leftmargin=*]   \setlength\itemsep{1em}
    \item  By definition of $v_{ij}(x)$, we have: $r v_{11}(x)-\mathcal{L}_{11} v_{11}(x) - f_{11}(x)= -rc_{12} >0,$ therefore $v_{ij}(x)$ satisfies (A21). On the other hand, we have:   $v_{11}(x)- \big (v_{21}(x) -c_{12} \big )=0,$ therefore $v_{ij}(x)$ satisfies (A22). With regards to (A23), we have:  $v_{11}(x)- \big (v_{12}(x) + \chi_{12} \big )=-\chi_{12} -\chi_{21} < 0,$ hence $v_{ij}(x)$ satisfies (A23). By Lemma \ref{lem:lemma3}, we deduce that $v_{ij}(x)$ solves (\ref{eq1}).
    \item By definition of $v_{ij}(x)$, we have: $v_{12}(x)- \big (v_{11}(x) + \chi_{21} \big )=0,$ hence $v_{ij}(x)$ satisfies (A31). On the other hand, we have:  $v_{12}(x)- \big (v_{22}(x) -c_{12} \big )=0,$  therefore $v_{ij}(x)$ satisfies (A32). By Lemma \ref{lem:lemma3}, we deduce that $v_{ij}(x)$ solves (\ref{eq2}).
    \item By definition of $v_{ij}(x)$, we have:  
    $r v_{11}(x)-\mathcal{L}_{11} v_{11}(x) - f_{11}(x)=0,$ hence $v_{ij}(x)$ satisfies (A11). On the other hand, we have: $v_{21}(x)- \big (v_{11}(x) -c_{21} \big )=c_{21}+c_{12} >0,$  hence $v_{ij}(x)$ satisfies (A12). With regards to (A13), we have: $v_{21}(x) - \big (v_{22}(x) +\chi_{12} \big )=-\chi_{12} -\chi_{21} < 0 ,$ therefore $v_{ij}(x)$ satisfies (A13). By Lemma \ref{lem:lemma3}, we deduce that $v_{ij}(x)$ solves (\ref{eq3}).
    \item By definition of $v_{ij}(x)$, we have:  $v_{22}(x)- \big (v_{21}(x) + \chi_{21} \big )=0,$  therefore $v_{ij}(x)$ satisfies (A32). On the other hand, $v_{22}(x)- \big (v_{21}(x) - c_{21} \big ) = -c_{21}-c_{12} < 0,$ therefore $v_{ij}(x)$ satisfies (A31). By Lemma \ref{lem:lemma3}, we deduce that $v_{ij}(x)$ solves (\ref{eq4}). \hfill $\Box$
\end{itemize}
\item  Case of \textbf{Condition B3:} 
The proof is similar as in \textbf{4.3 - 2.} 
\item Case of \textbf{Condition B4:} 
The proof is similar as in \textbf{4.3 - 4.}
 \end{enumerate}
\epf
\bethe
\normalfont
\label{th:theo4.5}
 Case where $K_{11} =   K_{21} >  K_{12} = K_{22}.$
\begin{enumerate}[leftmargin=*]
    \item If Condition B1 is verified then let $v_{ij}(x)$ be defined as follows:

\begin{tcolorbox}[colback=yellow!10!white, colframe=red!50!black]
\begin{minipage}{0.48\textwidth}
\begin{fleqn}[\parindent]
\begin{spreadlines}{1ex}
     \begin{align}
    v_{11}(x)& =\begin{cases}
    \hat{V}_{11}(x) +  Ax^{m^+_{11}}, &   x < x^*\\
\hat{V}_{12}(x)  +  \chi_{12},  &  x \geq x^*
\end{cases},\\
v_{12}(x)& = \hat{V}_{12}(x),\qquad \qquad \qquad x >0, 
     \end{align}
     \end{spreadlines}
     \end{fleqn}
\end{minipage}
\hfill\vline\hfill
\begin{minipage}{0.5\textwidth}
\begin{fleqn}[\parindent]
\begin{spreadlines}{1ex}
   \begin{align}
   v_{21}(x)& = \begin{cases}
\hat{V}_{21}(x) +  Ax^{m^+_{21}}, & x <x^*\\
\hat{V}_{12}(x)  +  \chi_{12},  & x \geq x^* 
\end{cases},\\
v_{22}(x)& =\hat{V}_{22}(x), \qquad \qquad \qquad x >0,
     \end{align}
     \end{spreadlines}
     \end{fleqn}
\end{minipage}
\end{tcolorbox}

 where $x^*=\inf S^{\eta}_{12} \in (0, +\infty)$. Using the fact that $v$ is $\mathcal{C}^1$ on $\partial S^{\xi}_{12}$ we get 
 $$(K_{11}-K_{12})(x^*)^{\gamma}=\frac{m_{11}^{+} \chi_{12}}{m_{11}^{+}- \gamma}, \quad
A=(K_{12}-K_{11})\frac{\gamma}{m_{11}^{+}}(x^*)^{\gamma -m_{11}^{+}}.$$
Such $v_{ij}(x)$ are the solutions associated to the system of quasi -variational inequalities with
\begin{center}
$S^{\xi}_{11}=S^{\xi}_{22}=S^{\eta}_{22}=S^{\xi}_{12}=S^{\eta}_{12}=S^{\xi}_{21}=\emptyset, \hspace{1 mm}  S^{\eta}_{11}= S^{\eta}_{21}=[x^*,+\infty)$
\end{center}
When the switching costs are positive (condition B1), player II has to stay in regime $1$ for x less than the threshold value $x^*$ and switch beyond that value. Player I on his side hass to never switch.
\item If Condition B2 is verified, then  let $v_{ij}(x)$ be defined as follows:
\begin{adjustwidth}{-1 cm}{-1 cm}
\begin{tcolorbox}[colback=yellow!10!white, colframe=red!50!black]
\begin{minipage}{0.48\textwidth}
\begin{fleqn}[\parindent]
\begin{spreadlines}{1ex}
     \begin{align}
  v_{11}(x) &= \begin{cases}
    \hat{V}_{21}(x) +  Ax^{m^+_{21}}-c_{12}, & x < x^*\\
\hat{V}_{22}(x)  -  c_{12} + \chi_{12}, & x\geq x^*
\end{cases},\\
\; v_{12}(x)& =\hat{V}_{22}(x) -c_{12}, \qquad \qquad \qquad x >0,
     \end{align}
     \end{spreadlines}
     \end{fleqn}
\end{minipage}
\hfill\vline\hfill
\begin{minipage}{0.5\textwidth}
\begin{fleqn}[\parindent]
\begin{spreadlines}{1ex}
   \begin{align}
  \;  v_{21}(x)& = \begin{cases}
    \hat{V}_{21}(x) + Ax^{m^+_{21}}, & \qquad x < x^*\\
\hat{V}_{22}(x)  +  \chi_{12} & \qquad x\geq x^*
\end{cases},  \\
\; v_{22}(x)& = \hat{V}_{22}(x), \qquad  \qquad \qquad \qquad x >0,
     \end{align}
     \end{spreadlines}
     \end{fleqn}
\end{minipage}
\end{tcolorbox}
\end{adjustwidth}
 where $x^*=\inf S^{\xi}_{21} \in (0, +\infty)$. Using the fact that $v$ is $\mathcal{C}^1$ on $\partial S^{\xi}_{12}$ we get :
$$(K_{21}-K_{22})(x^*)^{\gamma}=\frac{m_{11}^{+} \chi_{12}}{m_{11}^{+}- \gamma}, \quad
A=(K_{22}-K_{21})\frac{\gamma}{m_{11}^{+}}(x^{*})^{\gamma -m_{11}^{+}}.$$
Such $v_{ij}(x)$ are the solutions associated to the system of quasi-variational inequalities with:
\begin{center}
$S^{\xi}_{11}= S^{\xi}_{12}=(0, +\infty),\quad 
 S^{\eta}_{11}=S^{\eta}_{21}=[x^*, +\infty), \quad S^{\eta}_{12}= S^{\xi}_{21}= S^{\xi}_{22}=S^{\eta}_{22}=\emptyset.$
 \end{center}
When the switching cost from regime $1$ to regime $2$ is negative for player I (condition B2), he has to always switch from regime $1$ to regime $2$, player II has to stay in regime $1$ for $x$ less than the threshold value $x^*$ and switch beyond that value.
\item If Condition B3 is verified then let $v_{ij}(x)$ be defined as follows:

\begin{tcolorbox}[colback=yellow!10!white, colframe=red!50!black]
\begin{minipage}{0.48\textwidth}
\begin{fleqn}[\parindent]
\begin{spreadlines}{1ex}
     \begin{align}
    v_{11}(x) &= \hat{V}_{12}(x)+\chi_{12},\\
    v_{12}(x)& =  \hat{V}_{12}(x), 
     \end{align}
     \end{spreadlines}
     \end{fleqn}
\end{minipage}
\hfill\vline\hfill
\begin{minipage}{0.5\textwidth}
\begin{fleqn}[\parindent]
\begin{spreadlines}{1ex}
   \begin{align}
   \; v_{21}(x)& = \hat{V}_{22}(x)+\chi_{12} , \\
\; v_{22}(x)& = \hat{V}_{22}(x).
     \end{align}
     \end{spreadlines}
     \end{fleqn}
\end{minipage}
\end{tcolorbox}

Such $v_{ij}(x)$ are the solutions associated to the system of quasi-variational inequalities with

\qquad \qquad \quad  $S^{\xi}_{11}=S^{\xi}_{22}= S^{\eta}_{22}= S^{\xi}_{21}=S^{\xi}_{12}= S^{\eta}_{12}=\emptyset, \quad
 S^{\eta}_{11}=S^{\eta}_{21}=(0, +\infty).$
 
 When the switching cost from regime $1$ to regime $2$ is negative for player II (condition B3), he has to always switch from regime $1$ to regime $2$, Player I has to never switch.
  \item If Condition B4 is verified then let $v_{ij}(x)$ be defined as follows: 

\begin{tcolorbox}[colback=yellow!10!white, colframe=red!50!black]
\begin{minipage}{0.48\textwidth}
\begin{fleqn}[\parindent]
\begin{spreadlines}{1ex}
     \begin{align}
     v_{11}(x)& = \hat{V}_{22}(x) + \chi_{12} -c_{12},\\
 v_{12}(x)& =  \hat{V}_{12}(x) -c_{12}, 
     \end{align}
     \end{spreadlines}
     \end{fleqn}
\end{minipage}
\hfill\vline\hfill
\begin{minipage}{0.5\textwidth}
\begin{fleqn}[\parindent]
\begin{spreadlines}{1ex}
   \begin{align}
  \; v_{21}(x)& = \hat{V}_{22}(x) + \chi_{12},\\
\; v_{22}(x)& = \hat{V}_{22}(x).
     \end{align}
     \end{spreadlines}
     \end{fleqn}
\end{minipage}
\end{tcolorbox}

Such $v_{ij}(x)$ are the solutions associated to the system of quasi -variational inequalities with
\begin{center}
$S^{\xi}_{11}=S^{\eta}_{11}=S^{\xi}_{12}= S^{\eta}_{21}=(0, +\infty) \quad S^{\eta}_{12}=S^{\xi}_{21}=S^{\xi}_{22}= S^{\eta}_{22}=\emptyset.$
    \end{center}
When the switching cost from regime $1$ to regime $2$ is negative for both players (condition B4), they has to always switch from regime $1$ to regime $2$, 
\end{enumerate}
    \eethe

\bpf
\normalfont
\begin{enumerate}[leftmargin=*]
    \item Consider the case of condition (B1):  Let $v_{ij}(x)$ be defined as in Theorem \ref{th:theo4.5}-1.  
let us first check solutions for $x <  x^{*} $:
\begin{itemize}[label=$\star$, leftmargin=*]   \setlength\itemsep{1em}
    \item  By definition of $v_{ij}(x)$, we have:  $r v_{11}(x)-\mathcal{L}_{11} v_{11}(x) - f_{11}(x)= 0,$ hence $v_{ij}(x)$ satisfies (A11). On the other hand,  we have: $v_{11}(x)-\big (v_{21}(x)- c_{12} \big)= c_{12} >0,$ therefore, $v_{ij}(x)$ satisfies (A12). With regards to (A13), let $p$ be defined as follows: $p(x)=\hat{V}_{11}(x)-\hat{V}_{12}(x)+ Ax^{m^+_{11}}-\chi_{12}= (K_{11} -K_{12}) x ^{\gamma}+ Ax^{m^+_{11}}-\chi_{12} .$
A straightforward calculus shows that $p'(x)$ is positive, this implies that $p(x)$ is increasing. Furthermore, $p(x^*)=0 $, hence $p(x)<0$, subsequently, $v_{ij}(x)$ satisfies (A13). By Lemma \ref{lem:lemma3}, we deduce that $v_{ij}(x)$ solves (\ref{eq1}).
    \item By definition of $v_{ij}(x)$, we have:  $r v_{12}(x)-\mathcal{L}_{12} v_{12}(x) - f_{12}(x)= 0,$ hence $v_{ij}(x)$ satisfies (A11). On the other hand, we have: $v_{12}(x)-\big (v_{22}(x)- c_{12} \big )= c_{12} >0,$ therefore, $v_{ij}(x)$ satisfies (A12). With regards to (A13), we have: $v_{12}(x)- \big (v_{11}(x) + \chi_{21} \big )=(K_{12} -K_{11}) x ^{\gamma} - Ax^{m^+_{11}} - \chi_{12}.$
Let $q$ be defined as follows: $q(x)=(K_{12} -K_{11}) x ^{\gamma} - Ax^{m^+_{11}} - \chi_{21}.$  A straightforward calculus shows that $q'(x)$ is negative, hence  $q$ is decreasing. $q(0)=- \chi_{21}$, subsequently $q$ is clearly negative, and hence $v_{ij}(x)$ satisfies (A13). By Lemma \ref{lem:lemma3}, we deduce that $v_{ij}(x)$ solves (\ref{eq2}).
\end{itemize}
We can notice that $v_{11}(x)=v_{21}(x)$ and $v_{12}(x)=v_{22}(x)$, so the proof for (\ref{eq1}) applies to (\ref{eq3}) and the proof for (\ref{eq2}) applies to (\ref{eq4}).\\
Let's now check solution for $x \geq  x^{*} $:
\begin{itemize}[label=$\star$, leftmargin=*]   \setlength\itemsep{1em}
    \item By definition of $v_{ij}(x)$, we have: $r v_{11}(x)-\mathcal{L}_{11} v_{11}(x) - f_{11}(x)= \frac{ K_{12} -K_{11}}{K_{11}} x^{\gamma} + r \chi_{12}= \big [\frac{-m^+_{11}}{K_{11}(m^+_{11}-\gamma ) } \chi_{12} + r \chi_{12} \big ].$
A straightforward calculus shows that the above quantity is negative and hence $v_{ij}(x)$ satisfies (A32). With regards to (A31), we have: $v_{11}(x)-\big (v_{12}(x) + \chi_{12} \big )=\hat{V}_{12}(x) + \chi_{12} - \big ( \hat{V}_{12}(x) + \chi_{12} \big )=0.$ By Lemma \ref{lem:lemma3}, we deduce that $v_{ij}(x)$ solves (\ref{eq1}).
\item By definition of $v_{ij}(x)$, we have:  
$r v_{12}(x)-\mathcal{L}_{12} v_{12}(x) - f_{12}(x)= 0,$
hence $v_{ij}(x)$ satisfies (A11). On the other hand,  we have: $v_{12}(x)-\big (v_{22}(x)- c_{12} \big )= c_{12} >0.$ therefore, $v_{ij}(x)$ satisfies (A12). With regards to (A13), we have: $v_{12}(x)- \big (v_{11}(x) + \chi_{21} \big )= - \chi_{21}-\chi_{12}.$ By Lemma \ref{lem:lemma3}, we deduce that $v_{ij}(x)$ solves (\ref{eq2}).
\end{itemize}
Again we have $v_{11}(x)=v_{21}(x)$ and $v_{12}(x)=v_{22}(x)$, so the proof for (\ref{eq1}) applies to (\ref{eq3}) and the proof for (\ref{eq2}) applies to (\ref{eq4}). \hfill $\Box$
\item Consider the case of condition (B2):  Let $v_{ij}(x)$ be defined as in Theorem \ref{th:theo4.5}-2.  Let us check solutions for $x < x^*$:

\begin{itemize}[label=$\star$, leftmargin=*]   \setlength\itemsep{1em}
    \item By definition of $v_{ij}(x)$, we have:    $r v_{11}(x)-\mathcal{L}_{11} v_{11}(x) - f_{11}(x)=  r c_{12} > 0,$ hence $v_{ij}(x)$ satisfies (A21). On the other hand, we have:  $v_{11}(x)-\big (v_{21}(x) -c_{12} \big )=0,$ therefore $v_{ij}(x)$ satisfies (A22). With regards to (A23), we have:
    $v_{11}(x)-\big (v_{12}(x)+\chi_{12} \big )=(K_{11}- K_{12})x^{\gamma} + Ax^{m^+_{11}} -\chi_{12}$. Let $h$ be defined as follows: $h(x)=(K_{11}- K_{12})x^{\gamma} + Ax^{m^+_{11}} -\chi_{12}.$ A straightforward calculus shows that $h$ is increasing; on the other hand, $h(x^*)=0$, hence  $v_{ij}(x)$ satisfies (A23). By Lemma \ref{lem:lemma3}, we deduce that $v_{ij}(x)$ solves (\ref{eq1}).
    \item By definition of $v_{ij}(x)$, we have:  $r v_{12}(x)-\mathcal{L}_{12} v_{12}(x) - f_{12}(x)=  r c_{12} > 0,$ hence $v_{ij}(x)$ satisfies (A21). On the other hand, $v_{12}(x)-\big (v_{22}(x) -c_{12} \big )=0,$ therefore, $v_{ij}(x)$ satisfies (A22). With regards to (A23), let $p$ be defined as follows:
    $p(x)=v_{12}(x)-\big (v_{11}(x)+\chi_{21} \big )=-(K_{11}- K_{12})x^{\gamma} - Ax^{m^+_{11}} -\chi_{12} - \chi_{21}.$ We can see that $p(x)=-h(x) -\chi_{12}-\chi_{21}$ where $h$ is function from above; $p$ is clearly decreasing; on the other hand, $p(0)=-\chi_{21}<0$. By Lemma \ref{lem:lemma3}, we deduce that $v_{ij}(x)$ solves (\ref{eq2}).
    \item By definition of $v_{ij}(x)$, we have:   
    $r v_{21}(x)-\mathcal{L}_{21} v_{21}(x) - f_{21}(x)=  0,$
     hence $v_{ij}(x)$ satisfies (A11). On the other hand we have: $v_{21}(x)-\big (v_{11}(x)- c_{21} \big )= c_{12} +c_{21}>0,$ hence  $v_{ij}(x)$ satisfies (A12). With regards to (A13), we have: $v_{21}(x)- \big (v_{22}(x) + \chi_{12} \big )=(K_{21} -K_{22}) x ^{\gamma} +Ax^{m^+_{21}} - \chi_{12}$. By Lemma \ref{lem:lemma3}, we deduce that $v_{ij}(x)$ solves (\ref{eq3}).
Let $q$ be defined as follows: $q(x)=(K_{21} -K_{22}) x ^{\gamma} + Ax^{m^+_{11}} - \chi_{12}.$ We can see that $q(x)=h(x)$ where $h$ is the function from above. Hence $v_{ij}(x)$ satisfies (A13). 
    \item By definition of $v_{ij}(x)$, we have:   
    $r v_{22}(x)-\mathcal{L}_{22} v_{22}(x) - f_{22}(x)=  0,$ therefore $v_{ij}(x)$ satisfies (A11). On the other hand, $v_{22}(x)-\big (v_{12}(x)-c_{21} \big )= c_{21} + c_{12} >0,$ hence $v_{ij}(x)$ satisfies (A12). With regards to (A13) let $z$ be defined as follows: $z(x)= v_{22}(x)-\big (v_{21}(x)-\chi_{21} \big )=-(K_{21} -K_{22}) x ^{\gamma} - Ax^{m^+_{21}} - \chi_{21},$ We can see that $z(x)=p(x)$, where $p$ is the function from above. Hence $v_{ij}(x)$ satisfies (A13). By Lemma \ref{lem:lemma3}, we deduce that $v_{ij}(x)$ solves (\ref{eq4}).
    
\end{itemize}

Let us check solution for $x \geq x^*$:

\begin{itemize}[label=$\star$, leftmargin=*]   \setlength\itemsep{1em}
    \item By definition of $v_{ij}(x)$, we have: 
    $v_{11}(x)- \big (v_{12}(x)+\chi_{12} \big )=0$,
    hence $v_{ij}(x)$ satisfies (A31). On the other hand, $v_{11}(x)-\big (v_{21}(x)+ c_{12} \big )=0,$  therefore $v_{ij}(x)$ satisfies (A32). By Lemma \ref{lem:lemma3}, we deduce that $v_{ij}(x)$ solves (\ref{eq1}).
    \item By definition of $v_{ij}(x)$, we have:  $r v_{12}(x)-\mathcal{L}_{12} v_{12}(x) - f_{12}(x)=  r c_{12} > 0,$ therefore $v_{ij}(x)$ satisfies (A21). On the other hand, we have:  $v_{12}(x)-\big (v_{22}(x) -c_{12} \big )=0,$  therefore $v_{ij}(x)$ satisfies (A22). With regards to (A23), we have :
   $v_{12}(x)-\big (v_{11}(x) -\chi_{21} \big )=-(\chi_{21} + \chi_{12})<0.$ By Lemma \ref{lem:lemma3}, we deduce that $v_{ij}(x)$ solves (\ref{eq2}).
   \item By definition of $v_{ij}(x)$, we have: $r v_{21}(x)-\mathcal{L}_{21} v_{21}(x) - f_{21}(x) <0 ,$ therefore $v_{ij}(x)$ satisfies (A32). On the other hand, we have:
$v_{21}(x)-\big (v_{22}(x) -\chi_{12} \big )=0,$ therefore $v_{ij}(x)$ satisfies (A31). By Lemma \ref{lem:lemma3}, we deduce that $v_{ij}(x)$ solves (\ref{eq3}).
\item By definition of $v_{ij}(x)$, we have:
$r v_{22}(x)-\mathcal{L}_{22} v_{22}(x) - f_{22}(x) =0.$
therefore $v_{ij}(x)$ satisfies (A11). On the other hand we have: $v_{22}(x)-\big (v_{12}(x)-c_{21} \big )=c_{21} + c_{12} > 0,$ therefore $v_{ij}(x)$ satisfies (A12). With regards to (A13), we have:
$v_{22}(x)-\big (v_{21}(x)+\chi_{21}\big)=-\chi_{21}-\chi_{12}<0,$ therefore $v_{ij}(x)$ satisfies (A13). By Lemma \ref{lem:lemma3}, we deduce that $v_{ij}(x)$ solves (\ref{eq4}). \hfill $\Box$
\end{itemize}
\item Consider the case of condition (B3): Let $v_{ij}(x)$ be defined as in Theorem \ref{th:theo4.5}-3
\begin{itemize}[label=$\star$, leftmargin=*]   \setlength\itemsep{1em}
\item By definition of $v_{ij}(x)$, we have: $r v_{11}(x)-\mathcal{L}_{11} v_{11}(x) - f_{11}(x)= \frac{K_{12}-K_{11}}{K_{11}}x^{\gamma}+r\chi_{12}<0,$ therefore $v_{ij}(x)$ satisfies (A32). On the other hand, we have:  $v_{11}(x)-\big (v_{12}(x)+\chi_{12} \big )=0,$ therefore $v_{ij}(x)$ satisfies (A31). By Lemma \ref{lem:lemma3}, we deduce that $v_{ij}(x)$ solves (\ref{eq1}). 
\item By definition of $v_{ij}(x)$, we have: 
$r v_{11}(x)-\mathcal{L}_{11} v_{11}(x) - f_{11}(x)=0,$
therefore $v_{ij}(x)$ satisfies (A11). On the other hand, $v_{12}(x)-\big (v_{22}(x)-c_{12} \big )=c_{12}>0,$ therefore $v_{ij}(x)$ satisfies  (A12). With regards to (A13), we have: $v_{12}(x)-\big(v_{11}(x)+\chi_{21} \big)=-(\chi_{21}+\chi_{12})<0,$ therefore $v_{ij}(x)$ satisfies  (A13). By Lemma \ref{lem:lemma3}, we deduce that $v_{ij}(x)$ solves (\ref{eq2}).
\end{itemize}
We can see that $v_{11}(x)=v_{21}(x)$ and $v_{12}(x)=v_{22}(x)$, so the proof for (\ref{eq1}) applies to (\ref{eq3}) and the proof for (\ref{eq2}) applies to (\ref{eq4}).\hfill $\Box$
\item Consider the case of condition (B4):  Let $v_{ij}(x)$ be defined as in Theorem \ref{th:theo4.5}-4
\begin{itemize}[label=$\star$, leftmargin=*]   \setlength\itemsep{1em}
    \item  By definition of $v_{ij}(x)$, we have:  $v_{11}(x)-\big (v_{21}(x) - c_{12} \big )= \hat{V}_{12}(x) - c_{12} + \chi_{12} - \big (\hat{V}_{12}(x)+  \chi_{12} - c_{12} \big )= 0,$  therefore $v_{ij}(x)$ satisfies (A32). On the other hand, we have:  $v_{11}(x)-\big (v_{12}(x) + \chi_{12} \big )= \hat{V}_{12}(x) - c_{12} + \chi_{12} - \big (\hat{V}_{12}(x) - c_{12} + \chi_{12} \big ) = 0,$  therefore $v_{ij}(x)$ satisfies (A31).  By Lemma \ref{lem:lemma3}, we deduce that $v_{ij}(x)$ solves (\ref{eq1}). 
    \item By definition of $v_{ij}(x)$, we have:  $r v_{12}(x)-\mathcal{L}_{12} v_{12}(x) - f_{12}(x)= -r c_{12} > 0,$ therefore $v_{ij}(x)$ satisfies (A21). On the other hand, we have: 
$v_{12}(x)- \big (v_{22}(x) - c_{12} \big )=\hat{V}_{12}(x) - c_{12}-\big (\hat{V}_{12}(x) - c_{12} \big )=0,$  therefore $v_{ij}(x)$ satisfies (A22). With regards to (A23), we have:  $v_{12}(x)- \big (v_{11}(x) + \chi_{21} \big )=\hat{V}_{12}(x) - c_{12}- \big (\hat{V}_{12}(x) - c_{12}+\chi_{12}+\chi_{21} \big )<0, $ therefore $v_{ij}(x)$ satisfies (A23).  By Lemma \ref{lem:lemma3}, we deduce that $v_{ij}(x)$ solves (\ref{eq2}). 
    \item By definition of $v_{ij}(x)$, we have: $r v_{21}(x)-\mathcal{L}_{21} v_{21}(x) - f_{21}(x)=\frac{ K_{12} -K_{21}}{K_{21}} x^{\gamma} + r \chi_{12},$
The above quantity is negative, therefore $v_{ij}(x)$ satisfies (A32). On the other hand, we have: $v_{21}(x)-\big (v_{22}(x)- \chi_{12} \big )=\hat{V}_{12}(x)+\chi_{12}-(\hat{V}_{12}(x)+\chi_{12})=0,$ therefore $v_{ij}(x)$ satisfies (A31). By Lemma \ref{lem:lemma3}, we deduce that $v_{ij}(x)$ solves (\ref{eq3}). 
    \item By definition of $v_{ij}(x)$, we have: 
    $r v_{22}(x)-\mathcal{L}_{22} v_{22}(x) - f_{22}(x)=0$, therefore $v_{ij}(x)$ satisfies (A11). On the other hand, we have: 
$v_{22}(x) - \big (v_{12}(x) - c_{21} \big )=\hat{V}_{12}(x)- \big (\hat{V}_{12}(x)-c_{21} \big )=c_{21},$ therefore $v_{ij}(x)$ satisfies (A12). With regards to (A13), we have:  $v_{22}(x) - \big (v_{21}(x) + \chi_{21} \big )=\hat{V}_{12}(x)-\big (\hat{V}_{12}(x)+\chi_{12}+\chi_{21} \big )=-(\chi_{12}+\chi_{21})<0,$ thus $v_{ij}(x)$ satisfies (A13). By Lemma \ref{lem:lemma3}, we deduce that $v_{ij}(x)$ solves (\ref{eq4}). \hfill $\Box$
\end{itemize}
\end{enumerate}    
\epf

\section{Numerical Procedure}
\setcounter{equation}{0}
\renewcommand{\theequation}{5.\arabic{equation}}
In this part, we suggest a numerical procedure to compute threshold value in case we know qualitative structures of switching regions without knowing explicitly the value of the threshold values. To showcase our approach, We restrict ourselves for the particular case the switching regions are as follow: 
$$S^{\eta}_{12}=[x_{12},+\infty),\quad S^{\eta}_{21}=[0, y_{21}], \quad S^{\xi}_{12}=[x'_{12} , +\infty), \quad S^{\xi}_{21}=\emptyset, \quad \text{with  }  0 < y_{21} < x'_{12} < x_{12} .$$
To this end, let's introduce the following notations: for $a, \quad b, \quad x >0$
\begin{equation}
\tau_{a}=\inf \{t>0: X^x_t=a\},  \quad \tau_{ab}=\inf \{t>\tau_{a}: X^x_t=b\}.
\end{equation}
Consider the following  expectations functionals:
\begin{equation}
    R_1(x, a) =\mathbb{E} \bigg [  e^{-r\tau_{a}} \bigg ] , \quad R_2(x, a, b )  =\mathbb{E} \bigg [  e^{-r\tau_{ab}} \bigg ], \quad R_3(x, a,b )  =\mathbb{E} \bigg [  e^{-r\tau_{a}} \mathbb{1}_{\tau_{a} < \tau_{b}} \bigg ], 
\end{equation}
\begin{equation}
F_1(a ) =\mathbb{E}\bigg[\integ{0}{\tau_{a}}e^{-rs}f(X^x_{t})dt\bigg ], \text{ }  F_2(a,b)  =\mathbb{E}\bigg[\integ{\tau_{a}}{\tau_{b}}e^{-rs}f(X^x_{t})dt\bigg ], \text{ } F_3(a,b )  =\mathbb{E}\bigg[\integ{0}{\tau_{a} \wedge \tau_{b}}e^{-rs}f(X^x_{t})dt\bigg ].
\end{equation}
The computation of these expectations functionals can be found in the appendix of \textcolor{blue}{\cite{PVT}}.
\begin{itemize}[leftmargin=*, itemsep=1.5 pt]
\item \textbf{\underline{Assume that $x \geq x_{12}$:}}
then the optimal strategy in the case of regime (1, 1) is that both players have to immediately switch to regime 2 paying $c_{12}$ and  $\chi_{12}$ and then find the optimal strategy for $J_{22}$.The optimal strategy in the case of regime (1, 2) is that  player I has to immediately switch to regime 2 paying $c_{12}$ and then find the optimal strategy for $J_{22}$.
The optimal strategy in the case of regime (2, 1) is that  player II has to immediately switch to regime 2 paying $\chi_{12}$ and then find the optimal strategy for $J_{22}$.The optimal strategy in the case of regime (2, 2) is to let the process diffuse until it hits  $y_{21}$ at $\tau_{y_{21}}$ , player I then switch to regime 1 paying $c_{21}$ then we let the process diffuse until it hits $x_{12}$ then player I switch back to regime 1 paying $c_{12}$ and the process repeats itself.
\begin{equation}
\begin{split}
J_{22}(x, y_{21} , x'_{12} ,  x_{12})=\mathbb{E} \bigg [\integ{0}{\tau_{y_{21}}}e^{-rt}f(X^x_{t})dt -
e^{-r\tau_{y_{21}}} c_{21} + \integ{\tau_{y_{21}}}{\tau_{ y_{21}x_{12}}}e^{-rt}f(X^x_{t})dt \\  +  \bigg (J_{22}(x_{12} , y_{21} , x'_{12} ,  x_{12})-c_{12} \bigg )e^{-r\tau_{ y_{21}x_{12}}} \bigg ]. 
\end{split}
\end{equation}
We can determine $J_{22}(x_{12} , y_{21} , x'_{12} ,  x_{12})$ by taking   $x= x_{12}$. We get :
$$ J_{22}(x_{12} , y_{21} , x'_{12} ,  x_{12})= \big [F_1(y_{21}) + R_1(y_{21})c_{21} + F_2(\tau_{y_{21}}, \tau_{ y_{21}x_{12}}) -R_2(\tau_{y_{21}}, \tau_{ y_{21}x_{12}})c_{12} \big ] \big [(1- R_2(x, y_{21}, x_{12} )) \big ]^{-1}.$$ Once we know $J_{22}$, we know all the other by the following:
\begin{align}
J_{11}(x,y_{21}, x'_{12}, x_{12})&=-c_{12} + \chi_{12} + J_{22}(x,y_{21}, x'_{12}, x_{12}),\\
J_{12}(x,y_{21}, x'_{12}, x_{12})&=-c_{12}  + J_{22}(x,y_{21}, x'_{12}, x_{12}),\\
J_{21}(x,y_{21}, x'_{12}, x_{12})&= \chi_{12} + J_{22}(x,y_{21}, x'_{12}, x_{12}).
\end{align}
 \item \textbf{\underline{Assume that $x'_{12}\leq  x< x_{12}$:}} then the optimal strategy in the case of regime (1, 1) is that player II has to immediately switch to regime 2 paying   $\chi_{12}$ and then find the optimal strategy for $J_{12}$.The optimal strategy in the case of regime (1, 2) is that we let the process diffuse until it hits $x_{12}$ at $\tau_{x_{12}}$ then player I  switches to regime 2 paying $c_{12}$ and then find the optimal strategy for $J_{22}$.
The optimal strategy in the case of regime (2, 1) is that  player II has to immediately switch to regime 2 paying $\chi_{12}$ and then find the optimal strategy for $J_{22}$.The optimal strategy in the case of regime (2, 2) is to let the process diffuse until it hits  $y_{21}$ at $\tau_{y_{21}}$ , player I then switch to regime 1 paying $c_{21}$ then we let the process diffuse until it hits $x_{12}$ then player I switch back to regime 1 paying $c_{12}$ and the process repeats itself. 
\begin{equation}
\begin{split}
J_{22}(x, y_{21} , x'_{12} ,  x_{12})=\mathbb{E} \bigg [\integ{0}{\tau_{y_{21}}}e^{-rt}f(X^x_{t})dt-e^{-r\tau_{y_{21}}} c_{21} + \integ{\tau_{y_{21}}}{\tau_{y_{21}x_{12}}}e^{-rt}f(X^x_{t})dt \\  \eqskip \qquad  +  \bigg (J_{22}(x_{12} , y_{21} , x'_{12} ,  x_{12})-c_{12} \bigg )e^{-r\tau_{y_{21}x_{12}}} \bigg ].  \end{split}
\end{equation} 
From above, we already know $J_{22}(x_{12} , y_{21} , x'_{12} ,  x_{12})$, hence we know $J_{22}(x , y_{21} , x'_{12} ,  x_{12})$ and the other by the following relationship:
\begin{align}
J_{11}(x,y_{21}, x'_{12}, x_{12})&= \chi_{12} + J_{12}(x,y_{21}, x'_{12}, x_{12}),\\
 J_{12}(x,y_{21}, x'_{12}, x_{12})&= \mathbb{E}\bigg[\integ{0}{\tau_{x_{12}}}e^{-rt}f(X^x_{t})dt-
e^{-r\tau_{x_{12}}} \bigg (c_{12}+ J_{22}(x_{12} , y_{21} , x'_{12} ,  x_{12}) \bigg ) \bigg ],\\
J_{21}(x,y_{21}, x'_{12}, x_{12})&= \chi_{12} + J_{22}(x,y_{21}, x'_{12}, x_{12}).
\end{align}
    \item \textbf{ \underline{Assume that $y_{21} \leq  x< x'_{12}$:}} then the optimal strategy in the case of regime (1, 1) is that player I has to immediately switch to regime 2 paying   $\chi_{12}$ and then find the optimal strategy for $J_{12}$.The optimal strategy in the case of regime (1, 2) is that we let the process diffuse until it hits $x_{12}$ at $\tau_{x_{12}}$ then player I  switches to regime 2 paying $c_{12}$ and then find the optimal strategy for $J_{22}$.
The optimal strategy in the case of regime (2, 1) is that  player II has to immediately switch to regime 2 paying $\chi_{12}$ and then find the optimal strategy for $J_{22}$.The optimal strategy in the case of regime (2, 2) is to let the process diffuse until it hits  $y_{21}$ at $\tau_{y_{21}}$ , player I then switch to regime 1 paying $c_{21}$ then we let the process diffuse until it hits $x_{12}$ then player I switch back to regime 1 paying $c_{12}$ and the process repeats itself.
\begin{equation}
\begin{split}
 J_{22}(x,y'_{21} , y_{21} , x'_{12} ,  x_{12})=\mathbb{E} \bigg [\integ{0}{\tau_{y_{21}}}e^{-rt}f(X^x_{t})dt-
e^{-r\tau_{y_{21}}} c_{21} + \integ{\tau_{y_{21}}}{\tau_{y_{21}x_{12}}}e^{-rt}f(X^x_{t})dt
\\  + \bigg (J_{22}(x_{12},y'_{21} , y_{21} , x'_{12} ,  x_{12})-c_{12} \bigg )e^{-r\tau_{y_{21}}} \bigg ].
\end{split}
\end{equation} 
From above, we already know $J_{22}(x_{12} , y_{21} , x'_{12} ,  x_{12})$, hence we know $J_{22}(x , y_{21} , x'_{12} ,  x_{12})$ and the other by the following relationship:\\
\begin{equation}
J_{11}(x,y_{21}, x'_{12}, x_{12})=\mathbb{E} \bigg[\integ{0}{\tau_{x'_{12}}}e^{-rt}f(X^x_{t})dt+
e^{-r\tau_{x'_{12}}} \bigg (\chi_{12} + J_{12}(x'_{12} , y_{21} , x'_{12} ,  x_{12}) \bigg ) \bigg ],
\end{equation}
\begin{equation}
J_{12}(x,y_{21}, x'_{12}, x_{12})=\mathbb{E} \bigg [\integ{0}{\tau_{x_{12}}}e^{-rt}f(X^x_{t})dt-
e^{-r\tau_{x_{12}}} \bigg ( c_{12} + J_{22}(x_{12} , y_{21} , x'_{12} ,  x_{12}) \bigg ) \bigg ], 
\end{equation}
\begin{equation}
\begin{split}
    J_{21}(x,y_{21}, x_{12} ^{'}, x_{12}) = \mathbb{E} \bigg[\integ{0}{\tau_{x'_{12}} \wedge \tau_{y_{21}} }e^{-rt}f(X^x_{t})dt  
+e^{-r\tau_{x'_{12}}} \bigg (\chi_{12} + J_{22}(x'_{12} , y_{21} , x'_{12} ,  x_{12}) \bigg ) \mathbb{1}_{\tau_{x'_{12}} < \tau_{y_{21}}} \qquad \quad  \\
 + e^{-r\tau_{y_{21}}} \bigg (c_{21} + J_{11}(x'_{12} , y_{21} , x'_{12} ,  x_{12}) \bigg )\mathbb{1}_{ \tau_{y_{21}}\tau_{y_{21}}} < \tau_{x'_{12}}  \bigg ]. \qquad \qquad \qquad
 \end{split}
\end{equation}
\item \textbf{\underline{Assume that $x \leq y_{21}$:}} then the optimal strategy in the case of regime (1, 1) is that player I has to immediately switch to regime 2 paying   $\chi_{12}$ and then find the optimal strategy for $J_{12}$.The optimal strategy in the case of regime (1, 2) is that we let the process diffuse until it hits $x_{12}$ at $\tau_{x_{12}}$ then player I  switches to regime 2 paying $c_{12}$ and then find the optimal strategy for $J_{22}$.
The optimal strategy in the case of regime (2, 1) is that  player II has to immediately switch to regime 2 paying $\chi_{12}$ and then find the optimal strategy for $J_{22}$.The optimal strategy in the case of regime (2, 2) is to let the process diffuse until it hits  $y_{21}$ at $\tau_{y_{21}}$ , player I then switch to regime 1 paying $c_{21}$ then we let the process diffuse until it hits $x_{12}$ then player I switch back to regime 1 paying $c_{12}$ and the process repeats itself.
\begin{align}
J_{11}(x,y_{21}, x'_{12}, x_{12})=\mathbb{E}\big[\integ{0}{\tau_{x'_{12}}}e^{-rt}f(X^x_{t})dt+
e^{-r\tau_{x'_{12}}} \bigg (\chi_{12} + J_{12}(x'_{12},y_{21}, x'_{12}, x_{12}) \bigg ) \bigg ], \\
J_{12}(x,y_{21}, x'_{12}, x_{12})=\mathbb{E}\bigg[\integ{0}{\tau_{x_{12}}}e^{-rt}f(X^x_{t})dt+
e^{-r\tau_{x_{12}}} \bigg ( J_{22}(x_{12},y_{21}, x'_{12}, x_{12}) -c_{12} \bigg )\bigg],\\
J_{21}(x,y_{21}, x'_{12}, x_{12})= -c_{21}  + J_{11}(x,y_{21}, x'_{12}, x_{12}),\\
J_{22}(x, y_{21} , x'_{12} ,  x_{12})=-c_{21}  + J_{12}(x,y_{21}, x'_{12}, x_{12}).
\end{align}
We now find $v_{ij}(x)$ by solving :   
\begin{equation}
v_{ij}(x)=min_{x'_{12}} max_{(y_{21},x_{12})} J_{ij}(x,  y_{21} , x'_{12} ,  x_{12})=max_{(y_{21},x_{12})} min_{x'_{12}}J_{ij}(x,  y_{21} , x'_{12} ,  x_{12}).
\end{equation}
In case of only one player is switching, this situation degenerates into a standard optimization problem.
Optimal $(y_{21}^* ,  x_{12} ^{' *}, x_{12}^*)$ are then easly obtained such that:
$$(y_{21}^* ,  x_{12} ^{' *}, x_{12}^*)= (y_{21} ,   x_{12} ^{' *} ,  x_{12}),  0 < y_{21} < x_{12} ^{' *} < x_{12},    v_{ij}(x)= J_{ij}(x,  y_{21}^* , x_{12} ^{' *} ,  x_{12}^*).$$
\end{itemize}
\section{Graphic illustrations}
In this part, we illustrate the switching strategy for both players using the example of Theorem 4.3-2. For sake of readability, we mention that colors in the graphics aim at differentiating initial regime of each player. We use cyan when the initial regime of player I is regime 1, we use red when the initial regime of player I is regime 2;  we use violet (resp magenta) when the initial regime of player II is regime 1 (resp regime 2).
\\
\\
\\
\\
\\
\\
\\
\begin{figure}[H]
\caption{Graphic Illustration of Strategies in Theorem 4.3-ii \\}
\vspace{2 mm }
 \begin{tikzpicture}
 \draw[line width=1. 5 pt] (13,-4)--(0,-4)--(0,0)--(13,0);
 
 \node at (3.8,-2.2) {$x^*_A$} ;
 \node at (8.2,-2.2) {$x^*_B$};
 \node at (14.25,0) {$\color{cyan}R 1$};
 \node at (14.25,-1) {$ \color{red}R 2$};
  \node at (14.25,-3) {$ \color{violet}R 1$};
 \node at (14.25,-4) {$ \color{magenta} R 2$};
 \node at (-0.8,-0.7) {$P I$ \big \{ };
 \node at (-0.8,-3.5)  {$P II$ \big \{ };

 \node at (2,0.2) {switch};
 \node at (10.2,0.2) {continue};
 \node at (2,-2.8) {continue};
 \node at (.5,0.2) {continue};
  \node at (14.25,-2) {$x$ } ;
 \node at (5.5,-2.8) {continue};
 \node at (2,-4.2){continue} ;
  \node at (5.5,-4.2) {continue};
 \node at (2, -1.2) {continue};
  \node at (5.5,-1.2) { continue} ;
  \node at (10.2,-1.2) {switch };
  \node at (10.2,-2.8) {continue} ;
  \node at (10.2,-4.2){continue};
  
 \draw[very thick,<-] (13,-2)--(0,-2);
 \draw[line width=1.5 pt] (13,-1)--(0,-1);
  \draw[line width=1.5 pt] (13,-3)--(0,-3);
  \draw[ very thick, ->, violet] (1.5,0)--(1.5,-0.75);
  \draw[very thick, ->, magenta] (2.8,0)--(2.8,-0.75);
  \draw[ very thick, ->, violet] (9.5,-1)--(9.5,-0.25);
  \draw[very thick, ->, magenta] (10.8,-1)--(10.8,-0.25);
  \draw[thick,dashed] (3.5,-4)--(3.5,0);
  \draw[thick,dashed] (8,-4)--(8,0);
\end{tikzpicture}
   \begin{minipage}{0.5\textwidth}
     \includegraphics[width=8cm, height=5cm]{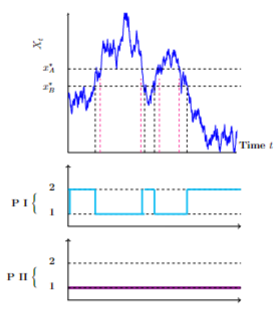}
     \caption{Theorem 4.3 B2: $I_0=(1,1)$ }\label{Fig:Data13}
   \end{minipage}
   \begin{minipage}{0.45\textwidth}
     \includegraphics[width=8cm, height=5cm]{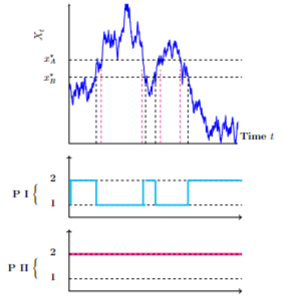}
     \caption{Theorem 4.3 B2: $I_0=(1,2)$}\label{Fig:Data14}
   \end{minipage}\\
   \begin{minipage}{0.5\textwidth}
     \includegraphics[width=8cm, height=5cm]{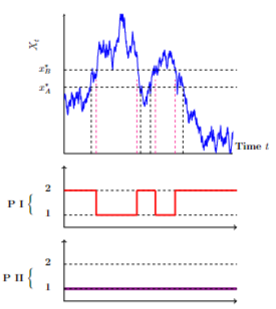}
     \caption{Theorem 4.3 B2: $I_0=(2,1)$ }\label{Fig:Data15}
   \end{minipage}
   \begin{minipage}{0.45\textwidth}
     \includegraphics[width=8cm, height=5cm]{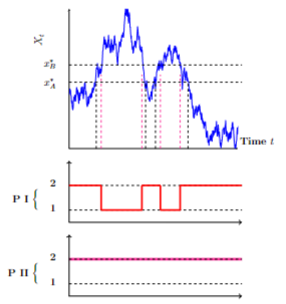}
     \caption{Theorem 4.3 B2: $I_0=(2,2)$}\label{Fig:Data16}
   \end{minipage}
\end{figure}

\no \textbf{Interpretation :} we can see from the above plots that, when cost of switching from regime 1 to regime 2 is negative for player I , since regime 1 pay equal reward, he has better stay in regime 1 except for lower values of $x$ where switching cost is more relevant than rewards themselves. He will gain in some way the opposite of the switching cost and hence greater than 0. Similarly, for larger values of $x$ where rewards are more relevant than switching costs,  he has to switch from regime 2 to regime 1. Player II on his side has to never switch since all regimes pay equal reward and switching costs are positives.

\newpage
\section{Conclusion}
In this paper, we have considered the problem of a two-player switching game characterized by a double obstacle HJBI equation. The first player's objective is to select a sequence of stopping times and a switching strategy so as to maximize the expected discounted cashflows, while the second player's objective is to select a sequence of stopping time and a switching strategy that minimizes the expected discounted cashflows. The aim was to derive an optimal explicit solution to the game in the assumption of identical profit functions and different diffusion operators. We characterize  the switching regions and turned the problem into one of finding a finite number of threshold values in the state process that trigger the switchings. Our major contribution is the derivation of explicit solution of the value function. We intend to develop this work by considering the case of different profit functions.\\ 

\no \textbf{Declarations}

\no \textbf{Conflict of interest.} The authors have not disclosed any competing interests.

\end{document}